\documentclass[12pt]{article} 
\long\def\kmcomment#1{}
\usepackage[leqno]{amsmath}
\usepackage{amsmath,etoolbox,mathtools}

\makeatletter
\newcommand{\leqnomode}{\tagsleft@true\let\veqno\@@leqno}
\newcommand{\reqnomode}{\tagsleft@false\let\veqno\@@eqno} 
\makeatother

\usepackage{setspace}
\usepackage{graphicx,psfrag}
\usepackage[svgnames]{xcolor}
\usepackage{tikz}
\usetikzlibrary{patterns,shapes}

\usepackage{newtxtext,newtxmath}
\usepackage[a4paper,truedimen,scale={.80,.85},centering]{geometry}

{\list{}{\leftmargin=#1 \topsep=0pt \parsep=0pt \itemsep=0pt}\item[]}%
{\endlist}
\RequirePackage{verbatim}
{\noindent\quote\endgraf\verbatim}%
{\endverbatim\endquote\endgraf\medskip}%
\makeatletter
\newenvironment{qverb*}%
{\noindent\quote\endgraf\@nameuse{verbatim*}}%
{\@nameuse{endverbatim*}\endquote\endgraf\medskip}%
\makeatother

\RequirePackage{fancyvrb} 
\DefineVerbatimEnvironment%
    {newverb}{Verbatim}
    {xleftmargin=10mm}

\DefineVerbatimEnvironment%
    {newverb*}{Verbatim}
    {xleftmargin=10mm,showspaces=true,showtabs=true} 

\usepackage{listings,quoting}
 \RequirePackage{listings} 
 \lstnewenvironment{kmverb}
    { \lstset{xleftmargin=10mm, breaklines=true}}
    { }
\newcommand{\kmqed}{\hfill\ensuremath{\blacksquare}}

\newcommand{\T}[1]{T_{#1}}

\newcommand{\myd}{d}

\newcommand{\futoji}[2]{{\mathbf #1}_{#2}}
\newcommand{\ds}{\ensuremath{\displaystyle }}
\newcommand{\pdel}{\partial}

\newcommand{\myCC}[2]{\overline{\text{C}}^{#1}_{#2}}

\newcommand{\myHH}[2]{\overline{\text{H}}^{#1}_{#2}}
 
\newcommand{\myHom}[1]{\text{H}_{#1}}
\newcommand{\myHomW}[2]{\text{H}_{#1,#2}}

\newcommand{\mR}{\ensuremath{\mathbb{R}}} 
\newcommand{\mZ}{\ensuremath{\mathbb{Z}}} 

\newcommand{\frakg}{\mathfrak{g}}

\newcommand{\tbdl}[1]{\mathrm{T}(#1)}

\newcommand{\Sbt}[2]{[#1,#2]}                
                
\newcommand{\mynabla}[2]{\nabla^{#1}_{#2}}

\usepackage{ntheorem} 
{
\theoremheaderfont{\bf}
\theorembodyfont{\rmfamily}
\theoremstyle{break}  
\theoremstyle{plain}  

\newtheorem{defn}{\bf Definition}
\newtheorem{prop}{Proposition}[section]
\newtheorem{exam}{Example}[section]
\newtheorem{remark}{Remark}[section]

\theoremstyle{plain}  

\newtheorem{theorem}{Theorem}[section] 
\newtheorem{thm}{Theorem}[section]
\newtheorem*{thm-none}{Theorem}[section]

\newtheorem{lemma}[theorem]{Lemma}
\newtheorem{kmCor}[theorem]{Corollary}

}

\renewcommand{\[}{$$} \renewcommand{\]}{$$}

\usepackage{newtxtext} 
\usepackage[enableskew,vcentermath]{youngtab}
\usepackage{listings,quoting}
\DefineShortVerb{\!} 

\newcommand{\myCS}[1]{ \text{C}_{#1}} 
\newcommand{\myOS}[1]{ \overline{\text{C}}_{#1}} 
\newcommand{\myCSW}[2]{ \text{C}_{#1,#2}} 
\newcommand{\myOSW}[2]{ \overline{\text{C}}_{#1,#2}} 
\newcommand{\mywedge}{\Delta} 
\AtBeginDocument{ \parindent=0pt }

\newcommand{\mathfrakX}[2]{\mathfrak{X}^{#1}_{#2}}
\newcommand{\bibun}[1]{(\frac{\pdel }{\pdel x_{#1}})} 
\newcommand{\AxA}[2]{A^{#1}_{#2}}
\newcommand{\XO}{E}

\numberwithin{equation}{section} 
\title{Euler number of homology groups of super Lie algebra}

\author{
Kentaro Mikami\thanks{ 
  Department of Computer Science and Engineering 
  Akita University, partially supported by
JSPS KAKENHI Grant Number  JP26400063, JP23540067 and JP20540059.}
 \and Tadayoshi Mizutani\thanks{
Professor Emeritus, Saitama University }
}

\allowdisplaybreaks[1]
\begin{document}
\maketitle



Main changes from the former version: 

(1) In the present file, we use the terminology ``pre Lie superalgebra'' or 
``\(\mZ\)-graded Lie superalgebra'' instead of ``super Lie algebra'' 
of the former version. 

(2) The title is changed: 
In the new file,  ``super Lie algebra'' is changed according  to (1) and 
added ``and Betti numbers''.  

(3)  Theorem 4.1 in the former version was not correct not 
including $0$-th chain. 
The revised one is Lemma \ref{thm:Euler:wZero} in this note. 

\kmcomment{
(4) New result (Theorem \ref{thm:m:th:Betti}): $m$-th Betti numbers
of $(w,h)$-double weighted chain complex is 0 if $w \ne h$.  

(5) New result (Theorem  \ref{thm:1:th:Betti}): The first Betti number 
of $(w,h)$-double weighted chain complex is 0.  
}

(4) Added new results:
For $(w,h)$-weighted chain complex, 
 $m$-th Betti numbers
 are 0 if $w \ne h$ (Theorem \ref{thm:m:th:Betti}) and  
the first Betti number is 0 
   (Theorem  \ref{thm:1:th:Betti}).

\bigskip

\section{Introduction} 
The well known de Rham cohomology group of a differentiable manifold $M$
is a cohomology group of the Lie algebra $\mathfrakX{ }{ }(M)$ of smooth
vector fields on $M$ together with the $\mathfrakX{ }{ }(M)$-module
\(\ds C^{\infty}(M)\) as coefficient.   Similarly, the Gel'fand-Fuks
cohomology theory is a cohomology theory of infinite dimensional Lie
algebras and there are many works on the cohomology of related subject
for example, the Lie algebra of volume preserving vector fields, the Lie
algebra of formal Hamiltonian vector fields and so on.  The notion of
these Lie algebra (co)homology groups is easy to understand, but the
calculation is hard to complete and one of the reason is  the infinity
of dimensions.  In order to reduce our computation to finite dimensional
case, we use an idea of ``weight'' (c.f.\,  for instance, \cite{M:N:K},
\cite{Mik:Nak}, \cite{KM:D6}, \cite{Kod:Mik:Miz:Torus}).   

There is (co)homology theory of  Lie superalgebras but few works of
\(\mZ\)-graded version.  Among Poisson geometers, \(\sum \Lambda^{p}
\tbdl{M} \) with the Schouten bracket is known as a prototype of
\(\mZ\)-graded (pre) Lie superalgebra and it is well-known that a
2-vector field \(\pi\) is Poisson if and only if \(\ds\Sbt{\pi}{\pi} =
0\).  Then Poisson condition \(\ds\Sbt{\pi}{\pi} = 0\) is equivalent to
\(\ds \pdel (\pi\wedge\pi) =0 \) in superalgebra homology theory, and
\(\ds\sqrt{ \ker (\pdel )}\) (the square root of cycles) is the space of
Poisson structures in some sense and there is some possibility of
studying Poisson structures in this direction. 

Thus, in this note, we will study homology groups of pre Lie
superalgebra  and relative homology groups with the coefficient
\(\frakg\)-module $V$ introducing (double) weight by many works of Lie
algebra (co)homology theory.  

\kmcomment{
studying  homology groups of pre Lie superalgebras
\(\frakg\) and relative homology groups with the coefficient
\(\frakg\)-module $V$ introducing (double) weight is very interesting.  

It is known the Schouten bracket plays a roll to characterize a 2-vector
field \(\pi\) is Poisson by \(\ds\Sbt{\pi}{\pi} = 0\).  Also, the graded
algebra of exterior power of tangent bundle is a pre Lie superalgebra
with the Schouten bracket.  
Then Poisson condition \(\ds\Sbt{\pi}{\pi} = 0\) is equivalent to \(\ds
\pdel (\pi\wedge\pi) =0 \).  Thus, \(\ds\sqrt{ \ker (\pdel )}\),   the
square root of cycles is the space of Poisson structures in some sense.
Thus, studying  homology groups of pre Lie superalgebras \(\frakg\) and
relative homology groups with the coefficient \(\frakg\)-module $V$
introducing (double) weight is very interesting.  }

\medskip

First we recall the definition of Lie superalgebra and pre Lie
superalgebra.  
\begin{defn}[(pre) Lie superalgebra]
Suppose a real vector space 
$\frakg $ is graded by \(\ds \mZ\) as 
\(\ds \frakg = \sum_{j\in \mZ} \frakg_{j} \)
and has a bilinear operation \(\Sbt{.}{.}\)  satisfying 
\begin{align}
& \Sbt{ \frakg_{i}}{ \frakg_{j}} \subset \frakg_{i+j} \label{cond:1} \\
& \Sbt{X}{Y} = (-1) ^{1+ x y} 
 \Sbt{Y}{X} \quad \text{ where }  X\in \frakg_{x} \text{ and } 
 Y\in \frakg_{y}  \\
& 
(-1)^{x z} \Sbt{ \Sbt{X}{Y}}{Z}  
+(-1)^{y x} \Sbt{ \Sbt{Y}{Z}}{X}  
+(-1)^{z y} \Sbt{ \Sbt{Z}{X}}{Y}  = 0 \quad \text{(Jacobi
identity).} \label{super:Jacobi}
\end{align}
Then we call \(\frakg\) a pre (or \(\mZ\)-graded) Lie superalgebra.   

A Lie superalgebra $\frakg $ is graded by \(\ds \mZ_{2}\) as \(\ds
\frakg = \frakg_{[0]} \oplus \frakg_{[1]} \) and the condition
\eqref{cond:1} is regarded as  \(\ds \Sbt{ \frakg_{[1]}}{ \frakg_{[1]}}
\subset \frakg_{[0]}\) in modulo 2 sense.   
\end{defn}

\begin{remark} 
Super Jacobi identity \eqref{super:Jacobi} above is equivalent to the one of
the following.  
\begin{align}
\Sbt{ \Sbt{X}{Y}}{Z} &=  
\Sbt{X}{ \Sbt{Y}{Z} } + (-1)^{y z} \Sbt{ \Sbt{X}{Z}}{Y} \\  
\Sbt{X}{\Sbt{Y}{Z}} &=  
\Sbt{\Sbt{X}{Y}}{Z}  + (-1)^{x y} \Sbt{Y}{ \Sbt{X}{Z}} 
\end{align} 

Suppose \(\ds \frakg = \sum_{j\in \mZ} \frakg_{j} \) is a pre Lie
superalgebra.  Let \(\ds\frakg_{[0]} = \sum_{i \text{ is even}}
\frakg_{i}\) and \(\ds\frakg_{[1]} = \sum_{i \text{ is odd}}
\frakg_{i}\). Then \(\ds \frakg = \frakg_{[0]} \oplus \frakg_{[1]} \)
holds and this is a Lie superalgebra.  
\end{remark}

\medskip

\begin{exam} \label{exam:super:elem} 
Take an $n$-dimensional vector space $V$ and split it as \(\ds V =
V_{0} \oplus V_{1}\).  Define \(\ds \frakg_{[i]} = \{ A\in
\mathfrak{gl}(V) \mid A (V_{j}) \subset V_{i+j} \} \). For each \(\ds A
\in \frakg_{[i]} \) and \(\ds B \in \frakg_{[j]} \), define
\(\ds\Sbt{A}{B} = A B - (-1)^{i j} B A \). Then  \(\ds \mathfrak{gl}(V)
=\frakg_{[0]} \oplus \frakg_{[1]} \) with this bracket is a Lie
superalgebra. 

More concretely, we take $n=2$ and \(\ds \dim V_{[i]}=1\) for \(i=0,1\).  
Then \(\ds \frakg_{[0]} = \begin{bmatrix} \ast & 0 \\ 0 & \ast
\end{bmatrix} \) and 
\(\ds \frakg_{[1]} = \begin{bmatrix}  0 & \ast  \\ \ast & 0
\end{bmatrix} \).    

Now define 
\(\ds \frakg_{0} =  \begin{bmatrix} a & 0 \\ 0 & -a \end{bmatrix}\), 
\(\ds \frakg_{1} = \begin{bmatrix}  0 & \ast  \\ \ast & 0
\end{bmatrix}\) and 
\(\ds \frakg_{2} =  \begin{bmatrix} a & 0 \\ 0 & a \end{bmatrix}\). Then 
\(\ds \mathfrak{gl}(2) =\frakg_{0}\oplus\frakg_{1}\oplus\frakg_{2}\) is
a pre Lie superalgebra.  
\end{exam}

\medskip

We will introduce the notion of double-weighted pre Lie superalgebras
(cf.\ Definition \ref{defn:w:weight}) and main results in this note are
the calculation of the Euler number of homology groups of
double-weighted pre  Lie superalgebras of special type.  
\kmcomment{
(cf.\ 
Lemma \ref{thm:Euler:wZero},  
Theorem \ref{thm:triv:general} and  
Theorem \ref{thm:module:gen}). 
}

\medskip
\begin{thm-none}(cf.\ 
Lemma \ref{thm:Euler:wZero})  
For general $n$, the Euler number of chain complex 
\(\ds \{ \myCSW{\bullet}{0,h}\} \) is $0$ for each $h$.  
\end{thm-none}

\begin{thm-none}(cf.\ 
Theorem \ref{thm:triv:general})  
For general $n$, the Euler number of chain complex 
\(\ds \{ \myCSW{\bullet}{w,h}\} \) 
is 0 for 
each $w$ and each $h$.  
\end{thm-none}


\begin{thm-none}(cf.\ Theorem \ref{thm:module:gen}) 
The Euler number of  
\(\ds (\myOSW{\bullet}{w,h},  \pdel_{V}) \) is 0 for each $w$ and $h$.   
\end{thm-none}

\begin{thm-none} (cf.\ Theorem \ref{thm:m:th:Betti})
The $m$-th Betti number of
\(\ds \{ \myCSW{\bullet}{w,h}\} \) 
is 0 for each double weight \((w,h)\) if \(w \ne h\). 
\end{thm-none}

\begin{thm-none} (cf.\ Theorem  \ref{thm:1:th:Betti})
The first Betti number of 
\(\ds \{ \myCSW{\bullet}{w,h}\} \) is 0 for each double weight \(w,h\). 
\end{thm-none}

\newcommand{\CY}[1]{ \widehat{\mathfrak Y}[#1]} %
\newcommand{\CYm}[2]{ \widehat{{\mathfrak Y}_{#2} }[#1]} %
\newcommand{\myL}[1]{\operatorname{L}_{#1}}

\section{Preliminaries, notations and basic facts} \label{sec:prelim} 
In the usual Lie algebra homology theory, $m$-th chain space is the
exterior algebra \(\ds \Lambda^{m}\frakg\) of \(\frakg\) and the
boundary operator is essentially comes from the operator \(\ds X \wedge
Y \mapsto \Sbt{X}{Y}\). 

In the case of pre Lie superalgebras, skew-symmetry of bracket operation
which yields "exterior algebra" is defined as the quotient of the tensor
algebra \(\ds \otimes^{m}  \frakg \) of \(\frakg \) by the two-sided
ideal generated by  \begin{align} & X \otimes Y + (-1)^{x y} Y \otimes X
\quad\text{where }\quad X \in \frakg_{x}, Y \in \frakg_{y} \;,
\end{align} and we denote the equivalence class of \(\ds X \otimes Y \)
by \(\ds X \mywedge Y\).  
\kmcomment{ \(\ds \bigtriangleup \triangle
\mywedge  \nabla \bigtriangledown \) } 
Since \( \ds X_{\text{odd}}
\mywedge Y_{\text{odd}}  = Y_{\text{odd}} \mywedge X_{\text{odd}} \) and
\( \ds X_{\text{even}} \mywedge Y_{\text{any}}  = - Y_{\text{any}}
\mywedge X_{\text{even}} \) hold, \( \mywedge ^{m} \frakg_{k} \) is a
symmetric algebra with respect to \(\mywedge\) for odd $k$ and is a skew-symmetric algebra 
with respect to \(\mywedge\) 
for even $k$.

\medskip
\begin{defn}
Assume that the pre Lie superalgebra \(\ds \frakg\) acts on a module
$V$ as follows: For each homogeneous \(\xi \in \frakg\) we have 
\(\ds \xi_{V} \in \text{End}(V) \) and satisfy 
\(\ds \Sbt{\xi}{\eta}_{V} = \xi_{V}  \circ \eta_{V}  
- (-1)^{|\xi| |\eta| } \eta_{V} \circ \xi_{V} \) where 
\(\ds \xi \in \frakg_{|\xi|}\) and 
\(\ds \eta \in \frakg_{|\eta|}\). We call $V$ be a \(\frakg\)-module.   
\kmcomment{
In short, 
we have a pre Lie superalgebra homomorphism  
\(\ds \frakg \ni \xi \mapsto \xi_{V} \in \text{End}(V)\).  
}
We often write \(\ds \xi_{V}(v) \) by
\(\ds \xi\cdot v\). 
\end{defn}
\begin{exam} A pre Lie superalgebra \(\ds \frakg\) is itself
\(\frakg\)-module by the own bracket
\(\ds X \cdot Z = \Sbt{X}{Z}\). 

Let $X,Y\in\frakg$ be homogeneous as \(\ds X\in \frakg_{x}, 
 Y\in \frakg_{y}\).    
\(\ds (X\circ  Y - (-1)^{x y} Y\circ X ) \cdot Z = 
\Sbt{X}{Y}\cdot Z \) holds and this is just Jacobi identity.  
\end{exam}

Suppose we have an exterior product of \(\ds Y_{1},\dots,Y_{m}\), i.e.,
\( Y_{1} \mywedge \cdots\mywedge Y_{m}\).  Omitting $i$-th element, we
have \(\ds Y_{1} \mywedge \cdots \mywedge Y_{i-1} \mywedge Y_{i+1}
\mywedge \cdots \mywedge Y_{m} \). It is often denoted as   \(\ds Y_{1}
\mywedge \cdots\widehat{Y_{i}}\cdots\mywedge Y_{m} \).  Here we denote
it by \(\ds \CYm{i}{m}  \).  If we omit $i$-th and $j$-th elements, we
denote the omitted product by \(\ds \CYm{i,j}{m}  \).  

\begin{defn} Let \(V\) be  a \(\frakg\)-module.   

For integer $m > 0$, define  
\(\ds \myOS{m} =\mywedge  ^{m}{\frakg} \otimes V = \sum_{i_{1}\leq
\ldots \leq i_{m}} \frakg_{i_{1}} \mywedge  \cdots \mywedge
\frakg_{i_{m}} \otimes V\), called $m$-th chain space.

In the case where $m=0$, we define \(\ds \myOS{0} = V\). 
\end{defn} 
We define a map \(\ds \pdel_{V} :\myOS{m}\to \myOS{m-1}\) by 
\begin{align} 
\pdel_{V} ( Y_{1}\mywedge \cdots \mywedge Y_{m}  \otimes v  ) = & 
\sum_{i<j} (-1)^{ 
{ \mathop{\sum}_{s<j} (1+ y_{j}y_{s})}  
+ { \mathop{\sum}_{s<i} (1+ y_{i}y_{s})}  } 
 \Sbt{Y_{j}}{Y_{i}} \mywedge \CYm{i,j}{m} \otimes v \label{pdelV:one}
  \\& 
+ (-1)^{m+1} 
\sum_{i=1}^{m} (-1)^{
{ \mathop{\sum}_{s>i} }
(1+ y_{i}y_{s})}  
 \CYm{i}{m} \otimes Y_{i} \cdot v \label{pdelV:two}
 \\\noalign{where \(\ds y_{i}\) is the degree of homogeneous element
 \(Y_{i}\), i.e., \(\ds Y_{i} \in \frakg_{y_{i}} \). } \notag
\end{align} 
We have the next basic fact.
\begin{thm}
\(\ds \pdel_{V}\circ \pdel_{V} = 0 \) holds. 
We have $m$-th homology group denoted by 
\[\ds 
\myHom{m}(\frakg, V) =  \ker(\pdel_{V} : \myOS{m} \rightarrow
\myOS{m-1})/ \pdel_{V} ( \myOS{m+1} )\;.  
\] 
\end{thm}

\begin{remark}
The first term of 
\eqref{pdelV:one} is also expressed as 
\begin{equation}
\sum_{i<j} (-1)^{ i-1 + y_{i} \mathop{\sum}_{i< s<j} y_{s}} 
Y_{1} \mywedge \cdots \widehat{ Y_{i} } \cdots \mywedge 
\underbrace{\Sbt{Y_{i}}{Y_{j}}}_{j} \mywedge \cdots  \mywedge Y_{m} \otimes v 
\;. 
\end{equation} 
At first, \(\ds Y_{i}\) moves to the left side of \(\ds Y_{j}\), then  the
parity changes to \(\ds (-1)^{ \sum_{s=i+1}^{j-1} ( 1 + y_{i}y_{s} )}
\). Then the ``side effect'' of the bracket operation produces \( \ds
(-1)^{j-1}\).  In this note, we have chosen  \( \ds (-1)^{j-1}\), but it
may be possible to choose \( \ds (-1)^{j}\). Then they have just the
opposite sign.  
\end{remark}

\begin{remark}
If \(\frakg\)-action on \(V\) is trivial, namely  \(\ds Y\cdot v=0 \)
for \(\ds \forall Y \in \frakg\) and \(\ds \forall v \in V\), then   
\eqref{pdelV:two} is always 0 and we may assume $V=\mR$. We
call this module the trivial module.  
Thus, when we essentially  deal with the trivial module, 
the chain space \(\ds \myCS{m} =
\Delta^{m} \frakg \) and \(\ds \pdel_{V}\) is \eqref{pdelV:one}
without $v$, which we denote \(\pdel\). It is clear that \(\ds
\pdel\circ \pdel = 0\) and we have the homology groups 
\[\ds 
\myHom{m}(\frakg, \mR) =  \ker(\pdel : \myCS{m} \rightarrow
\myCS{m-1})/ \pdel ( \myCS{m+1} )\;.  
\] 
\end{remark}

\subsection{Homology groups weighted by the first grading} 
Assume that a $\frakg$-module $V$ is \(\mZ\)-graded, i.e, \(\ds V =
\sum_{i} V_{i}\),  and satisfies \(\ds \frakg_{i}\cdot V_{j} \subset
V_{i+j} \). 
\begin{defn}
We define a non-zero element in 
\(\ds
\frakg_{i_{1}} \mywedge  \cdots \mywedge
\frakg_{i_{m}} \otimes V_{j}\) to have \(\ds i_{1}+\dots + i_{m} + j\)
as the (first) weight.  Define the subspace of 
\(\myOS{m} \) by 
\(\ds \myOSW{m}{w} = 
 \sum_{
\substack{
i_{1}\leq \ldots \leq i_{m}\\
\sum_{s=1}^{m}i_{s} +j = w
}
} 
\frakg_{i_{1}} \mywedge  \cdots \mywedge
\frakg_{i_{m}} \otimes V_{j}\), which is the direct sum of different
types of spaces of elements with weight $w$.  
\end{defn}

\begin{prop}
The (first) weight $w$ is preserved by \(\ds \pdel_{V}\), i.e., we have 
\(\ds \pdel_{V}( \myOSW{m}{w} ) \subset \myOSW{m-1}{w} \). Thus, we have
for a fixed $w$, 
$w$-weighted homology groups
\[\ds 
\myHomW{m}{w} (\frakg, V) =  \ker(\pdel_{V} : \myOSW{m}{w} \rightarrow
\myOSW{m-1}{w})/ \pdel ( \myOSW{m+1}{w} )\;.  
\] 
When $V$ is the trivial module, then we have 
\[\ds 
\myHomW{m}{w} (\frakg, \mR) =  \ker(\pdel_{} : \myCSW{m}{w} \rightarrow
\myCSW{m-1}{w})/ \pdel ( \myCSW{m+1}{w} )\;.  
\] 
\end{prop}
\subsection{Double-weighted homology groups} 
\begin{defn}[Double-weight] \label{defn:w:weight}
Assume that each subspace \(\ds \frakg_{i} \) of a given pre Lie superalgebra
\(\frakg\) is directly decomposed by subspaces \( \ds \frakg_{i,j}\)
as \(\ds \frakg_{i}  = \sum_{j} \frakg_{i,j} \) and satisfies 
\begin{equation}
\Sbt{ X }{ Y } \in \frakg_{i_1+i_2, j_1+j_2} \quad \text{for each}\ 
 X\in \frakg_{i_1,j_1} \ ,  \ 
Y\in \frakg_{i_2,j_2} \;. \end{equation}  
We say such pre Lie superalgebras are double-weighted. 

Assume that $\frakg$-module $V$ is also double-weighted \(\ds V_{i,j}\)
and satisfies \(\ds \frakg_{i,j} \cdot  V_{i',j'}
\subset V_{i+i', j+j'} \). 

Then we may define double-weighted $m$-th chain space by
\[
\myOSW{m}{w,h} = \sum_{
\substack{ 
i_{1}\leq \ldots \leq i_{m}\;,\; 
\sum_{s=1}^{m}i_{s} +i_{0} = w \\
\sum_{s=1}^{m}h_{s} +h_{0} = h 
}
} 
\frakg_{i_{1},h_{1}} \mywedge  \cdots \mywedge
\frakg_{i_{m},h_{m}} \otimes V_{i_{0},h_{0}} 
\]
\end{defn}

\begin{prop}
The double-weight $(w,h)$ is preserved by \(\ds \pdel_{V}\), i.e., we have 
\(\ds \pdel_{V}( \myOSW{m}{w,h} ) \subset \myOSW{m-1}{w,h} \). Thus, we have 
$(w,h)$-weighted homology groups
\[\ds 
\myHomW{m}{w,h} (\frakg, V) =  \ker(\pdel_{V} : \myOSW{m}{w,h} \rightarrow
\myOSW{m-1}{w,h})/ \pdel ( \myOSW{m+1}{w,h} )\;.  
\] 
When $V$ is the trivial module, then we have 
\[\ds 
\myHomW{m}{w,h} (\frakg, \mR) =  \ker(\pdel_{} : \myCSW{m}{w,h} \rightarrow
\myCSW{m-1}{w,h})/ \pdel ( \myCSW{m+1}{w,h} )\;.  
\] 
\end{prop}

\section{Pre Lie superalgebras with the Schouten bracket} 
A prototype of pre Lie superalgebra is the exterior algebra of 
the sections of exterior power of  tangent bundle of a differentiable manifold $M$ of
dimension $n$
\[ \frakg = \sum_{i=1}^{n} \Lambda^{i} \tbdl{M} = 
\sum_{i=0}^{n-1} \frakg_{i} 
\;, \quad \text{where}\quad 
\frakg_{i} = \Lambda^{i+1} \tbdl{M} \;  
\] with the Schouten bracket. 

There are several ways of defining 
the Schouten bracket, namely, axiomatic explanation, sophisticated one
using Clifford algebra or more direct ones (cf.\ 
\cite{Mik:Miz:homogPoisson}).  Here in the context of Lie algebra homology
theory, we introduce the Schouten bracket as follows:

\begin{defn}[Schouten bracket]
For \(\ds A \in \Lambda^{a} \tbdl{M} \) and 
\(\ds B \in \Lambda^{b} \tbdl{M} \), we define a binary operation
\(\Sbt{\cdot}{\cdot}\) by 
\begin{equation}\label{defn:abst:schouten}
(-1)^{a+1} \Sbt{A}{B} = 
 \pdel (A\wedge B) - (\pdel A)\wedge B -
(-1)^{a} A \wedge \pdel B \;. 
\end{equation}
\end{defn}

In some sense, the Schouten bracket measures how far from the derivation
the boundary operator \(\pdel\) is. 

The first chain space is
\(\ds \myCS{1} = \frakg= \sum_{p=1}^{n} \Lambda^{p} \tbdl{M}\).  
The second chain space is 
\begin{align*} \myCS{2} = \frakg \mywedge  \frakg  
= \sum_{p \leq q } \Lambda^{p} \tbdl{M} \mywedge \Lambda^{q} \tbdl{M}  
=&  \Lambda^{1} \tbdl{M} \mywedge \Lambda^{1} \tbdl{M} + 
  \Lambda^{1} \tbdl{M} \mywedge \Lambda^{2} \tbdl{M} + \cdots 
\\& 
 + \Lambda^{2} \tbdl{M} \mywedge \Lambda^{2} \tbdl{M}
 + 
  \Lambda^{2} \tbdl{M} \mywedge \Lambda^{3} \tbdl{M} + \cdots 
\end{align*}  
\begin{remark} \label{remark:Poisson}
Let \(\ds \pi\in\Lambda^{2}\tbdl{M}\). Then \(\ds \pi\mywedge\pi \in 
 \Lambda^{2} \tbdl{M} \mywedge \Lambda^{2} \tbdl{M} \subset \myCS{2}\)
 and \(\ds \pdel ( \pi\mywedge\pi ) = \Sbt{\pi}{\pi} \in \myCS{1}\).
 Thus, 
\(\ds \pi\in\Lambda^{2}\tbdl{M}\) is Poisson if and only if 
 \(\ds \pdel ( \pi\mywedge\pi ) = 0 \), and we express it by \(\ds \pi
 \in \sqrt{ \ker(\pdel) }\) symbolically. 
It will be interesting to study   
\(\ds \sqrt{ \ker(\pdel) }\) and also interesting to study specific
properties of Poisson structures in 
\(\ds \sqrt{ \pdel( \myCS{3} ) } \), which come from the boundary image of 
the third chain space  \(\myCS{3}\).  
\end{remark}
\kmcomment{
In particular, \(\ds 
\frakg_{-1} = \Lambda^{0} \tbdl{M} = C^{\infty}(M) \). 
}


In this pre Lie superalgebra, possible weights are non-negative integers.  
When weight is 0,  the chain spaces with trivial action 
are simply given by 
\(\ds \myCSW{m}{0} = \mywedge ^{m} \frakg_{0}= 
\mywedge ^{m} \tbdl{M} 
\)
and the homology is the Lie algebra homology of vector fields 
for \(m=1,\dots,n \).
For lower weights 1 or 2, 
 the chain spaces are simply given by 
\begin{align*}
\myCSW{m}{1} & = \mywedge^{m-1} \frakg_{0} \mywedge  \frakg_{1} = 
\mywedge^{m-1} \tbdl{M} \mywedge  
\Lambda^{2} \tbdl{M}\quad \text{for}\quad m=1,\dots  \;,
\\
\myCSW{m}{2} & = \mywedge^{m-1} \frakg_{0} \mywedge  \frakg_{2} 
\oplus \mywedge^{m-2} \frakg_{0} \mywedge^{2}\frakg_{1} 
\\& 
= 
\mywedge^{m-1} \tbdl{M} \mywedge  \Lambda^{3} \tbdl{M}
\oplus \mywedge^{m-2} \tbdl{M} \mywedge^{2}\Lambda^{2} \tbdl{M} 
\quad \text{for}\quad m=1,\dots  \;.
\end{align*} 
\begin{remark}
In particular, 
\(\ds \myCSW{1}{2}  = \Lambda^{3} \tbdl{M}\), 
\(\ds \myCSW{2}{2}  =
\tbdl{M} \mywedge  \Lambda^{3} \tbdl{M}
\oplus \Lambda^{2} \tbdl{M} \mywedge\Lambda^{2} \tbdl{M} \),  
\(\ds \myCSW{3}{2}  =
\tbdl{M} \mywedge  \tbdl{M} \mywedge  
\Lambda^{3} \tbdl{M}
\oplus \tbdl{M} \mywedge\Lambda^{2} \tbdl{M}  \mywedge\Lambda^{2} \tbdl{M} 
\).  Thus, by introducing weight, the chain spaces become smaller and
research becomes a little clear and easier.   
\end{remark}

Given a general weight $w$, the sequences 
\( 0\leq i_{1}\leq \dots\leq i_{m}\leq n-1\) with 
\( \sum_{s=1}^{m}  i_{s} = w\) 
correspond to \(\ds 1\leq j_{1}\leq \dots \leq j_{m}\leq n\) with 
\(\sum_{s=1}^{m} j_{s} = m+ w\) 
by \( j_{s} = 1+i_{s}\).  It is known that the  non-increasing sequences \(\ds
j_{m},\dots,j_{1}\) are Young diagrams of area \(w+m\) and length \(m\)
and we get the  original sequences by \(\ds i_{s} = j_{s} -1 \).  
\kmcomment{
When we pay attention to multiplicity \(\ds j_{i}\) of 
\(\ds \mywedge^{j_{i}} \Lambda^i \tbdl{M}\) 
then \( k_{i}\) is given by 
\(\ds k_{i} = \# \{s \mid j_{s}= i\} \) 
for each Young diagram \(\ds (j_{m},\ldots,j_{1}) \). 
}
\kmcomment{

(1) For $0\leq  i_ 1\leq \cdots \leq i_m \leq n-1$, 
 putting $ j_s=i_s+1$ we obtain an obvious sequence 
$j_m \geq \cdots \geq j_1 $.
This sequence gives a Young diagram of length $m$ and area 
$w+m$.  
(2)}
For each Young diagram $\{j_m,\dots,j_1\}$,
looking at the 'multiplicity'$j_i $ in 
$\Delta^{j_i}\Lambda^i T(M) $, we obtain a sequence 
$[k_1,k_2,\dots,k_n] $ consisting of 
$ k_i= \#\{ s \vert j_s=i\}$. 
\kmcomment{
(2')
The 2nd expression of the diagram $\{j_m,\dots,j_1\}$ 
is the sequence $[k_1,k_2,\dots,k_n] $ where 
$k_i$ is  the number of 
the 'multiplicity'$j_i $ of 
$\Delta^{j_i}\Lambda^i T(M) $.
Here, $n $ is equal to $j_m $.
}

\kmcomment{
of \(\ds \lambda \)  with length $m$ and area $w+m$. 
}

\begin{remark}[3 ways of Young diagram]
A Young diagram \(\lambda\) is a non-decreasing sequence of positive
integers, say \(\ds a_{1} , \dots, a_{m}\). For instance,
\(\tiny \yng(4,1,1)\) is a sequence of \(4,1,1\), here we denote it
as
\( {}^t (4,1,1) \) where superscript $t$ means ``traditional expression''.  
As explained above, when we focus to multiplicity of elements, we have
another sequence, in the concrete example above, \(2,0,0,1\) and we denote
it by \([2,0,0,1]\). Sometimes we have to write many 0 in this expression.  
The 3rd expression of Young diagram is measuring the height of each
column from left  to right. Again in the concrete example, we have a
sequence \(3,1,1,1\) and denote it by \( \langle  3,1,1,1\rangle \) and call it tower
(vertical) decomposition.  It is known in
general that the sequence of tower decomposition of \(\lambda\)
is just the conjugate of  \(\lambda \), i.e,  
\(\langle \lambda\rangle = {}^t( \text{conjugate of } \lambda) \).  
In detail of relations of those, refer to \cite{Mik:Miz:homogPoisson}.  
\end{remark}

\begin{remark}
We remark that $m$ does not stop at \(\dim M\)
in general because of property of our new ``wedge product'' \(\triangle\).  
\end{remark}

\section{Euler number of homology groups of concrete pre Lie superalgebras}\label{sec:concrete}
In the previous section, we have pre Lie superalgebras for each
differentiable manifold $M$.  In this section, 
we consider the Euclidean space \(\ds M = \mR^
{n} \) with the Cartesian coordinates \(x_{1},\dots,x_{n}\). 
Then, we get a pre Lie super subalgebra consisting of multi vector fields of
polynomial coefficients. 
We define \begin{align*} \frakg_{i,j} & =
\mathfrakX{i+1}{j+1} (\mR^{n}) 
 = \{
(i+1)\text{-multi vector fields with }
(j+1)\text{-homogeneous polynomials}\} \;. \end{align*} 
We see easily that \(\ds \Sbt{ \frakg_{i_1,j_1}} 
{ \frakg_{i_2,j_2}} \subset \frakg_{i_1+i_2, j_1+j_2} \) and so we get a
double-weighted pre Lie superalgebra.  
The spaces 
\(\ds \frakg_{i,j} \) are finite dimensional, precisely 
\(\ds \dim \frakg_{i,j} = \tbinom{n-1+j+1}{n-1} \tbinom{n}{i+1} \), and  
we study each component of  \(\ds \myCSW{m}{w,h} \) 
in the next subsection. 

\newcommand{\SubCC}[2]{\operatorname{SubC}^{(#1)}({#2})}
\newcommand{\SubCD}[2]{\operatorname{SubC}^{(#1)}_{[#2]}}

\subsection{Double weighted chain space \(\ds \myCSW{m}{w,h}\) }

\begin{prop}\label{prop:chainspace:general}
The chain space 
\(\ds \myCSW{m}{w,h} 
 = \sum_{ } 
\mathfrakX{i_{1}}{h_{1}} (\mR^{n}) \mywedge \cdots \mywedge 
\mathfrakX{i_{m}}{h_{m}} (\mR^{n}) \) 
of the double-weighted pre Lie superalgebra above is characterized as
follows:

\begin{enumerate}
\item
\(\ds (i_{s})_{s=1}^{m} \) are non-descending 
sequences of sum $w+m$ and length $m$. Since each entry
is less than $n+1$, we may count the multiplicity of them as follows:  
\([ k_{1},\dots,k_{n} ]\) where \(\ds k_{a} = \# \{ s \mid i_{s} = a
\}\) or denote it by $a\!:\!k_{a}$.  
\[\ds (i_{1}, \dots, i_{m}) = ( 
 \underbrace{1,\dots,1}_{k_{1}}, 
\dots,
 \underbrace{n,\dots,n}_{k_{n}} 
 ) = 
(1\!:\!{ k_{1}},\dots,n\!:\!{k_{n}} ) = [ k_{1},\dots, k_{n}] 
 \; , \]
we have 
\[
\sum_{s=1}^{n} k_{s} = m\;,\; \text{and}\; \sum_{s=1}^{n} s k_{s} = w + m \;.  
\]
Now denote 
\(\ds 
\mathfrakX{i_{1}}{h_{1}} (\mR^{n}) \mywedge \cdots \mywedge 
\mathfrakX{i_{m}}{h_{q}} (\mR^{n}) \) by \(\ds  
\mathfrakX{( i_{1} \leq \dots \leq i_{m})}{
( h_{1},\dots, h_{m})}
= \mathfrakX{[ k_{1}, \dots ,k_{n}]}{
( h_{1},\dots, h_{m})}\). 

\item
Each \(\ds h_{s}\) is non-negative integer and 
\( \sum_{s=1}^{m} (h_{s}-1) = h  \), and so 
\( \sum_{s=1}^{m} (h_{s}+1) = h +2 m  \).  
Consider the sequences of Young diagrams of area $w+2 m$ and length $m$.
Since those \(\{ h_{s}+1\}\) are not necessarily non-decreasing, we need to get all permutations of them,
then shift 1 negatively simultaneously.  
\item
 Assume \(\ds i_{p-1} < i_{p} = \dots = i_{q} =k < i_{q+1}\).  
Then we may relabel so that \(\ds h_{p} \leq \dots \leq h_{q} \), 
we write 
\[
\mathfrakX{k}{h_{p}} (\mR^{n}) \mywedge \cdots \mywedge 
\mathfrakX{k}{h_{q}} (\mR^{n})  = 
 \SubCC{k:(q-p+1)}{h_{p},\ldots,h_{q}} \]

\item
  Assume \(\ds i_{p} = \dots = i_{q} \) and \(\ds h_{p} = \dots =
h_{q} \). Then 
\[\ds 
\mathfrakX{i_{p}}{h_{p}} (\mR^{n}) \mywedge \cdots \mywedge 
\mathfrakX{i_{q}}{h_{q}} (\mR^{n})  = 
 \SubCC{i_{p}:(q-p+1)}{
 \underbrace{ h_{p},\ldots,h_{p} }_{q-p+1}} = 
\mywedge^{q-p+1}  \mathfrakX{i_{p}}{h_{p}} (\mR^{n}) \]
\begin{itemize}
\item
If \(\ds i_{p}\) is even, then  
\(\ds 
\mywedge^{q-p}  \mathfrakX{i_{p}}{h_{p}} (\mR^{n}) \) is a symmetric
algebra and its dimension is 
\[ \binom{ \tbinom{n}{i_p}  \tbinom{n-1+h_p}{n-1}-1+q-p+1 } {q-p+1} \]

\item
If \(\ds i_{p}\) is odd, then  
\(\ds 
\mywedge^{q-p+1}  \mathfrakX{i_{p}}{h_{p}} (\mR^{n}) \) is a skew-symmetric
algebra and its dimension is
\[ \binom{ \tbinom{n}{i_p}  \tbinom{n-1+h_p}{n-1} } {q-p+1} \]

In particular, if \(\ds q-p+1 > 
 \tbinom{n}{i_p}  \tbinom{n-1+h_p}{n-1}  
\) then the algebra is
0-dimensional. 
\end{itemize}
\end{enumerate}
\end{prop}

Introducing a new notation 
\begin{equation} \label{eqn:new:notation}
\SubCD{ k:\ell} {u} =   
\mathop{ \oplus}_{ \sum_{s=1}^{\ell}( h_{s}-1)= u} 
 \SubCC{ k:\ell} {h_{1},\dots,h_{s} } \;, \end{equation} 
and using the notations in Proposition \ref{prop:chainspace:general}, 
we have 
\begin{kmCor}\label{kmcor:one} 
\begin{equation}
 \myCSW{m}{w,h} = \sum_{
\substack{ 
 \sum_{i=1}^{n} k_{i} = m \\ 
 \sum_{i=1}^{n} i k_{i} = w+m\\ 
 \sum_{i=1}^{n}  u_{i} = h
}} 
\SubCD{ 1:k_{1}} {u_{1}} \mywedge\cdots\mywedge 
\SubCD{ n:k_{n}} {u_{n}}
\label{eqn:myCSW:SubCD}
\end{equation} 

\end{kmCor}
\kmcomment{
Define \(\ds \kappa_{i} = \dim \mathfrakX{1}{i-1} (\mR^{n}) \) for
$i=1,\dots$. 
}

\begin{prop} \label{prop:test0}
Assume $k$ is an odd integer.  
Let \([ \ell_{1}, \ell_{2},\dots ]\) be the sequence of multiplicities
of \(\ds h_{1}+1,\dots, h_{m}+1 \), where \(\ds \ell_{b} = \# \{ i \mid
h_{i} + 1 = b \} \).   
Then 
\[
\ds\SubCC{k:m}{h_1,\dots,h_m} = \mywedge ^{\ell_{1}} \mathfrakX{k}{0} 
 \mywedge ^{\ell_{2}} \mathfrakX{k}{1} 
 \mywedge \cdots 
\] holds.  If an inequality  
\(\ds  \ell_{i} \leq \dim   \mathfrakX{k}{i-1}  = \tbinom{n}{k}
\tbinom{n-1+i-1}{n-1} \) holds for each $i$, then 
\(\ds\SubCC{k:m}{h_1,\dots,h_m} \) is non trivial, and whose dimension is  
\(\ds
\prod _{i} \binom{ 
 \tbinom{n}{k} \tbinom{n-1+i-1}{n-1}}{ \ell_{i}}  
\). 
\end{prop}

\textbf{Proof:} Since $k$ is odd, each algebra 
\(\ds \mywedge ^{\ell_{i}} \mathfrakX{k}{i-1} \) is skew-symmetric and
 the proposition holds comparing the dimension of  
\(\ds \mathfrakX{k}{i-1} \). \kmqed 


The requirements of the chain space 
\(\ds \myCSW{m}{w,h} \) 
in \eqref {eqn:myCSW:SubCD} are  
$w=\sum_{s=1}^{m} (s-1)k_{s} $ and 
 $\sum_{s=1}^{m} k_{s} =m $
for  the first weight.  Thus, 
if $w=0$, then \(\ds [k_{1},k_{2},\dots]= [m,0,\dots]\) or 
if $w=1$, then \(\ds [k_{1},k_{2},\dots]= [m-1,1,0,\dots]\).  
If $w=2$, then \(\ds [k_{1},k_{2},\dots]= [m-2,2,0,\dots]
\) or \( [m-1,0,1, 0, \dots]\). 
  
From the definition of the second weight $h$, we see that \(\ds
\sum_{s=1}^ m h_s = m + h\), and so \( m \geq -h\), in more precise, 
\( m \geq \max(-h,1)\). About upper bound of the range $m$, we discuss
later. 


\subsubsection{The first weight \(w=0\) case }
In this subsection, we assume $w=0$. Then the algebra is just Lie
algebra and  we see that 
\begin{align*} 
\myCSW{m}{0,h} & = 
\SubCD{ 1:m} {h}  = \sum_{ \substack{ \sum_{t} \ell_{t} = m \\
 \sum_{t} t \ell_{t} = 2m + h }}
\mywedge ^{\ell_{1}} \mathfrakX{1}{0} 
 \mywedge ^{\ell_{2}} \mathfrakX{1}{1} 
 \mywedge \cdots 
\kmcomment{
 \quad\text{
where}\quad [\ell_{1},\ell_{2},\dots ] = {}^{t} (
h_{1},\ldots,h_{m})\;.}
\;. 
\end{align*} 
We have some restrictions from the proposition \ref{prop:test0}. 

\begin{prop} Assume $w=0$ and \(\ds \myCSW{m}{0,h} \ne 0 \).  

Then  $h \geq  - \dim \mathfrakX{1}{0} $, and also  
\(
\max(1,-h) \leq m \leq h+ 2 \dim \mathfrakX{1}{0} + 
\dim \mathfrakX{1}{1} \) holds. 
\end{prop}
\textbf{Proof:} We follow the notation above, then 
\begin{align}
\sum_{t} \ell_{t} &= m \label{eq:length}  \\
 \sum_{t} t \ell_{t} &= 2m + h \label{eq:area}
\end{align}
Since \( \eqref{eq:area}  - 2 \eqref{eq:length}\), we have 
\(\ds -\ell_{1} + \sum_{s>2}(s-2) \ell_{s} = h \), 
and 
\(\ds \sum_{s>2}(s-2) \ell_{s} = h+ \ell_{1}\), thus we have \( 0 \leq
h+ \ell_{1}\).  Applying 
the requirement \(\ds \ell_{1} \leq \dim \mathfrakX{1}{0}\), we have  
\( 0 \leq h+  \dim \mathfrakX{1}{0}\).

From \( \eqref{eq:area}  - 3 \eqref{eq:length}\), we have 
\( -2\ell_{1} -\ell_{2} + \sum_{s>3}(s-3) \ell_{s} = -m + h \), thus 
\( 0 \leq \sum_{s>3}(s-3) \ell_{s} = -m + h + 2 \ell_{1}
+ \ell_{2} \). 
Applying the requirement \(\ds \ell_{2} \leq \dim \mathfrakX{1}{1}\), we have  
\( m-h \leq 2 \dim \mathfrakX{1}{0} + \dim \mathfrakX{1}{1} 
\). 
\kmqed

\begin{remark} 
In the previous proposition, we have an upper bound of $m$. If we use
the third or more higher comparison, we have more sharp estimate of
upper bound of $m$.  
\end{remark}

\begin{exam} \label{exam:n2:w0}
Assume $n=2$ for simplicity and we study the chain
space 
\begin{align*} 
\myCSW{m}{0,h} &  
= \sum_{ \substack{ \sum_{t} \ell_{t} = m \\
 \sum_{t} t \ell_{t} = 2m +h }}
\mywedge ^{\ell_{1}} \mathfrakX{1}{0} 
 \mywedge ^{\ell_{2}} \mathfrakX{1}{1} 
 \mywedge \cdots 
 \quad .
\end{align*} 
Assume $h = -2$. Then $m$ starts from 2.  The possible Young diagrams
are characterized by area $2m-2$ and length $m$.  We see that the Young
diagram \( \langle  m,m-2 \rangle = [ 1^2, 2^{m-2} ]\) is only candidate
for our chain space.  
Thus 
\(\ds \myCSW{m}{0,-2} = 
\mywedge ^{2} \mathfrakX {1}{0}  (\mR^{2}) 
\mywedge ^{m-2} \mathfrakX {1}{1}  (\mR^{2}) \) 
and we get dimension for each space as follows: the Euler number is 0. 
\begin{center}
\begin{tabular}{c | *{6}{c} }
m & 2 & 3 & 4 & 5 & 6 \\\hline
$\dim$ & 1 & 4 & 6 & 4 & 1 \\
$\dim\pdel$ & 1 & 3 & 3 & 1 & 0 \\\hline
Betti & 0 &  0 &  0 &  0 &  0  
\end{tabular}
\hfil 
\(
\begin{array}{c | *{6}{c} }
 &       m-1  & m & m+1 \\\hline
\dim &  d_{m-1} & d_{m} & d_{m+1} \\
\dim\pdel  & r_{m-1} & r_{m} & r_{m+1} \\\hline
\text{Betti} &    & d_{m} -( r_{m-1}+r_{m})   & \\\hline   
\dim\ker\pdel  & k_{m-1} & k_{m} & k_{m+1} \\
\text{Betti} & & ( k_{m}+k_{m+1})- d_{m+1} & \\   
\end{array}
\)
\end{center}

Assume $h = -1$.  The area is $2m-1$, and the good Young diagrams are $
\langle m,m-1\rangle $ or $\langle  m,m-2 , 1\rangle  $ and so \(\ds [
1^{1}, 2^{m-1} ]\) or \(\ds [ 1^2, 2^{m-3}, 3^1 ]\).  
Thus  
\[\ds \myCSW{m}{0,-1} = 
\mathfrakX {1}{0}  (\mR^{2}) 
\mywedge ^{m-1} \mathfrakX {1}{1}  (\mR^{2}) 
\oplus 
\mywedge^{2}
\mathfrakX {1}{0}  (\mR^{2}) 
\mywedge ^{m-3} \mathfrakX {1}{1}  (\mR^{2}) 
\mywedge \mathfrakX {1}{2}  (\mR^{2}) 
\;.
\] 
So we get dimension for each space as follows: the Euler number is 0. 
\begin{center}
\begin{tabular}{c | *{7}{c} }
m & 1& 2 & 3 & 4 & 5 & 6 & 7 \\\hline
$\dim$ & 2 &  8 & 18 & 32 & 38 & 24  & 6 \\
$\dim\pdel$ & 2 &  6 & 12 & 20 & 18 & 6  & 0 \\\hline
Betti & 0 &  0 & 0 & 0 & 0 & 0  & 0 
\end{tabular}
\end{center}

Assume $h = 0$.  The area is $2m$ and good Young diagrams are $\langle
m,m\rangle  $, $\langle  m,m-1, 1\rangle $,  $\langle  m,m-2, 1, 1\rangle$ 
or $\langle  m , m-2 , 2\rangle$ 
and so \(\ds [1^0, 2^{m} ]\),  
 \(\ds [ 1^1, 2^{m-2}, 3^1]\), 
 \(\ds [ 1^2, 2^{m-3}, 4^1]\) or  
 \(\ds [ 1^2, 2^{m-4}, 3^2]\).   
Thus 
\begin{align*}\ds \myCSW{m}{0,0} & = 
\mywedge^{m}
\mathfrakX {1}{1}  (\mR^{2}) 
\oplus \mathfrakX {1}{0}  (\mR^{2}) 
\mywedge ^{m-2} \mathfrakX {1}{1}  (\mR^{2}) 
\mywedge \mathfrakX {1}{2}  (\mR^{2}) 
\\&\quad 
\oplus 
\mywedge^{2} \mathfrakX {1}{0}  (\mR^{2}) 
\mywedge^{m-3} \mathfrakX {1}{1}  (\mR^{2}) 
\mywedge \mathfrakX {1}{3}  (\mR^{2}) 
\oplus 
\mywedge^{2} \mathfrakX {1}{0}  (\mR^{2}) 
\mywedge^{m-4} \mathfrakX {1}{1}  (\mR^{2}) 
\mywedge^{2} \mathfrakX {1}{2}  (\mR^{2}) \; .
\end{align*} 
When $h=0$, zero-th chain space is defined and  
\(\ds \myCSW{0}{0,0}  = \mR\). Thus,  
the dimension for each space and the rank of \(\pdel\) are as follows: the Euler number is $0$. 
\begin{center}
\begin{tabular}{c | *{9}{c} }
m & 0 & 1& 2 & 3 & 4 & 5 & 6 & 7 & 8 \\\hline
$\dim$ &1& 4 &  18 & 60 & 120 & 156 & 134  & 68 & 15 
\\
$\dim\pdel$ & 0 & 4 & 14 & 46 & 74 & 80 & 54 & 13 & 0 
\\\hline
Betti & 1 &  0 & 0 & 0 & 0 & 2 & 0 & 1 & 2 
\end{tabular}
\end{center}

Assume $h = 1$. The possible Young diagrams have area $2m+1$ and length
$m$ and decompositions are $\langle  m,m,1\rangle $, $\langle m,
m-1,1,1\rangle $, $\langle  m,m-1 , 2\rangle $, $\langle
m,m-2,1,1,1\rangle$, $\langle m,m-2,2,1\rangle$,  $\langle
m,m-2,3\rangle$ 
and so \(\ds [1^{0}, 2^{m-1}, 3^1 ]\),  
 \(\ds [ 1^1, 2^{m-2},3^0, 4^1]\), 
 \(\ds [ 1^1, 2^{m-3}, 3^2]\), 
 \(\ds [ 1^2, 2^{m-3},3^0,4^0, 5^1]\) or  
 \(\ds [ 1^2, 2^{m-4}, 3^1, 4^1]\),  
 \(\ds [ 1^2, 2^{m-5}, 3^3]\). 
Thus 
\begin{align*}\ds \myCSW{m}{0,1} & = 
\mywedge^{0} \mathfrakX {1}{0}  (\mR^{2}) 
\mywedge^{m-1} \mathfrakX {1}{1}  (\mR^{2}) 
\mywedge^{} \mathfrakX {1}{2}  (\mR^{2}) 
\oplus \mathfrakX {1}{0}  (\mR^{2}) 
\mywedge ^{m-2} \mathfrakX {1}{1}  (\mR^{2}) 
\mywedge^{0} \mathfrakX {1}{2}  (\mR^{2}) 
\mywedge \mathfrakX {1}{3}  (\mR^{2}) 
\\&\quad 
\oplus \mathfrakX {1}{0}  (\mR^{2}) 
\mywedge ^{m-3} \mathfrakX {1}{1}  (\mR^{2}) 
\mywedge^{2} \mathfrakX {1}{2}  (\mR^{2}) 
\oplus \mywedge^{2} \mathfrakX {1}{0}  (\mR^{2}) 
\mywedge ^{m-3} \mathfrakX {1}{1}  (\mR^{2}) 
\mywedge \mathfrakX {1}{4}  (\mR^{2}) 
\\&\quad 
\oplus \mywedge^{2} \mathfrakX {1}{0}  (\mR^{2}) 
\mywedge ^{m-4} \mathfrakX {1}{1}  (\mR^{2}) 
\mywedge^{ } \mathfrakX {1}{2}  (\mR^{2}) 
\mywedge^{ } \mathfrakX {1}{3}  (\mR^{2}) 
\oplus \mywedge^{2} \mathfrakX {1}{0}  (\mR^{2}) 
\mywedge ^{m-5} \mathfrakX {1}{1}  (\mR^{2}) 
\mywedge^{3} \mathfrakX {1}{2}  (\mR^{2}) 
\end{align*}
The dimension for each space is as follows: the Euler number is $0$. 
\begin{center}
\begin{tabular}{c | *{9}{c} }
m & 1& 2 & 3 & 4 & 5 & 6 & 7 & 8 & 9 \\\hline
$\dim$ & 6 &  40 & 140 & 328 & 522 & 544  & 352 & 128 &  20 \\
$\dim\pdel$ & 6 &34 & 106 & 222 & 300 & 244 & 108 &20&0 \\\hline
Betti  & 0 &  0 & 0 & 0 & 0 & 0  & 0 & 0 &  0 
\end{tabular}
\end{center}

\kmcomment{
When $h = 2$, the good Young diagrams with 
 area $2m+2$ and length $m$ are 
$\langle  m,m,1,1\rangle $,  
$\langle  m,m,2\rangle $,  
$\langle  m,m-1,1,1,1\rangle $,  
$\langle  m,m-1,2,1\rangle $,  
$\langle  m,m-1,3\rangle $,  
$\langle  m,m-2,1,1,1,1\rangle $,  
$\langle  m,m-2,2,1,1\rangle $,  
$\langle  m,m-2,2,2\rangle $,  
$\langle  m,m-2,3,1\rangle $,  
$\langle  m,m-2,4\rangle $,  

And so the good diagrams are 
$[0,m-1,0,1]$,  
$[0,m-2,2]$,  
$[ 1,m-2,0,0,1 ]$,  
$[ 1,m-3,1,1 ]$,  
$[ 1,m-4,3] $,  
$[ 2,m-3,0,0,0,1] $,  
$[ 2,m-4,1,0,1] $,  
$[ 2,m-4,0,2] $,  
$[ 2,m-5,2,1] $,  
$[ 2,m-6,4] $.    
Thus 
\begin{align*}\ds \myCSW{m}{0,2} = &\;   
\mywedge^{0}\mathfrakX{1}{0}
\mywedge^{m-1}\mathfrakX{1}{1}\mywedge^{0}\mathfrakX{1}{2}\mywedge^{1}\mathfrakX{1}{3} 
+  \mywedge^{0}\mathfrakX{1}{0}
\mywedge^{m-2}\mathfrakX{1}{1}\mywedge^{2}\mathfrakX{1}{3} 
\\ & 
+  \mywedge^{ 1}\mathfrakX{1}{0}
\mywedge^{m-2}\mathfrakX{1}{1}\mywedge^{0}\mathfrakX{1}{2}\mywedge^{0}\mathfrakX{1}{3}\mywedge^{1}
\mathfrakX{1}{4} 
+  \mywedge^{ 1}\mathfrakX{1}{0}
\mywedge^{m-3}\mathfrakX{1}{1}\mywedge^{1}\mathfrakX{1}{2}\mywedge^{1}
\mathfrakX{1}{3} 
\\ & 
+  \mywedge^{1}\mathfrakX{1}{0}
\mywedge^{m-4}\mathfrakX{1}{1}\mywedge^{3}\mathfrakX{1}{2} 
+  \mywedge^{ 2}\mathfrakX{1}{0}
\mywedge^{m-3}\mathfrakX{1}{1}\mywedge^{0}\mathfrakX{1}{2}\mywedge^{0}\mathfrakX{1}{3}\mywedge^{0}\mathfrakX{1}{4}\mywedge^{1}\mathfrakX{1}{5} 
\\ & 
+  \mywedge^{ 2}\mathfrakX{1}{0}
\mywedge^{m-4}\mathfrakX{1}{1}\mywedge^{1}\mathfrakX{1}{2}\mywedge^{0}\mathfrakX{1}{3}\mywedge^{1}\mathfrakX{1}{4} 
+  \mywedge^{2}\mathfrakX{1}{0}
\mywedge^{m-4}\mathfrakX{1}{1}\mywedge^{0}\mathfrakX{1}{2}\mywedge^{2}\mathfrakX{1}{3} 
\\ & 
+  \mywedge^{2}\mathfrakX{1}{0}
\mywedge^{m-5}\mathfrakX{1}{1}\mywedge^{2}\mathfrakX{1}{2}\mywedge^{1}\mathfrakX{1}{3} 
+  \mywedge^{ 2}\mathfrakX{1}{0}
\mywedge^{m-6}\mathfrakX{1}{1}\mywedge^{4}\mathfrakX{1}{2} 
\end{align*}

The dimension for each space is as follows: the Euler number is $0$. 
\begin{center}
\begin{tabular}{c | *{10}{c} }
m & 1& 2 & 3 & 4 & 5 & 6 & 7 & 8 & 9 & 10 \\\hline
$\dim$ & 8 & 67 & 296 & 802 & 1428 & 1730  & 1400 & 698 &  180 & 15 \end{tabular}
\end{center} 
}
\end{exam}

We discussed in \cite{Mik:Miz:homogPoisson} the Euler number of Lie
algebra homology groups of given weight $w$ and homogeneity $h$ of
Poisson tensor where we dealt Young diagrams of area $w+(2-h)m$ with
length $m$.  By discussion there, we have next result.  
\begin{lemma}\label{thm:Euler:wZero}
For general $n$, the Euler number of chain complex 
\(\ds \{ \myCSW{\bullet}{0,h}\} \) is $0$. 
\end{lemma}
\textbf{Proof:} We use the notations in \cite{Mik:Miz:homogPoisson}.  
Since 
\begin{align*} 
\myCSW{m}{0,h} &  
= \sum_{ \substack{ \sum_{t} \ell_{t} = m \\
 \sum_{t} t \ell_{t} = 2m + h }}
\mywedge ^{\ell_{1}} \mathfrakX{1}{0} 
 \mywedge ^{\ell_{2}} \mathfrakX{1}{1} 
 \mywedge \cdots \quad ,
\kmcomment{
 \text{
where}\quad [\ell_{1},\ell_{2},\dots ] = {}^{t} ( h_{1},\ldots,h_{m})\;.
}
\end{align*} 
we have to deal with Young diagrams \(\ds \mynabla{2m+h}{m}\)  of area
$2m+h$ with length $m$.   A recursive formula 
\begin{equation}
\label{eqn:recursive}
\mynabla{2m+h}{m} = B\cdot \mynabla{2m+h-1}{m-1} \sqcup
\T{m} \mynabla{m+h}{m}  
\end{equation}
holds.
If $h=0$ then we have 
\[\ds 
\mynabla{2m}{m} = B\cdot \mynabla{2m-1}{m-1} \sqcup
\T{m} \mynabla{m+h}{m}  
= \T{m}\cdot \T{m}  \sqcup B\cdot \mynabla{2m-1}{m-1} \;.
\]  
Thus, \(\ds \dim \myCSW{m}{0,0} = \tbinom{\dim \mathfrakX{1}{1}}{m} 
+ \sum_{ \lambda \in \mynabla{2m-1}{m-1}} \dim (B\cdot \lambda) \).  
 
When we denote each 
\( \lambda \in \mynabla{2m-1}{m-1} \) by \(\ds [ \ell_{1},\ell_{2},
\dots]\),  \( \ds \sum_{t}\ell_{t} = m-1 \) and  
 \( \ds \sum_{t} t \ell_{t} = 2m-1 \) must be satisfied. 
 \(\ds B\cdot \lambda =  [1+ \ell_{1},\ell_{2},\dots ] \) and so   
 \(\ds \dim (B\cdot \lambda) = 
 \tbinom{\dim \mathfrakX{1}{0}}{1+ \ell_{1}}  
\tbinom{\dim \mathfrakX{1}{1}}{ \ell_{2}} 
\tbinom{\dim \mathfrakX{1}{2}}{ \ell_{3}}\dots  
 \).  

We easily see that \(\ds \sum_{m>0} (-1)^{m}  \tbinom{\dim
\mathfrakX{1}{1}}{m} = -1\). 
About the alternating sum of the second term, 
we have 
\begin{align*} 
& 
\sum_{m>0} (-1)^{m} \sum_{ 
\lambda \in \mynabla{2m-1}{m-1} 
} \dim (B\cdot \lambda) \\ =  & 
\sum_{m>0} (-1)^{m} \sum_{ 
\substack{
\sum_{t}\ell_{t} = m-1 \\ 
\sum_{t} t \ell_{t} = 2m-1} 
}
 \tbinom{\dim \mathfrakX{1}{0}}{1+ \ell_{1}}  
\tbinom{\dim \mathfrakX{1}{1}}{ \ell_{2}} 
\tbinom{\dim \mathfrakX{1}{2}}{ \ell_{3}}\dots  
\\ =& 
\sum (-1)^{ 1+ \sum_{t}\ell_{t}  } 
\sum_{ 
2(1+\sum_{t}\ell_{t})  = 1+ \sum_{t} t \ell_{t} 
}
 \tbinom{\dim \mathfrakX{1}{0}}{1+ \ell_{1}}  
\tbinom{\dim \mathfrakX{1}{1}}{ \ell_{2}} 
\tbinom{\dim \mathfrakX{1}{2}}{ \ell_{3}}\dots  
\\ =& 
\sum (-1)^{ \ell_{2}  } 
\tbinom{\dim \mathfrakX{1}{1}}{ \ell_{2}} 
\sum_{ 
2(1+\sum_{t}\ell_{t})  = 1+ \sum_{t} t \ell_{t} 
} (-1)^{ 1+ \sum_{t \ne 2}\ell_{t}  } 
 \tbinom{\dim \mathfrakX{1}{0}}{1+ \ell_{1}}  
\tbinom{\dim \mathfrakX{1}{2}}{ \ell_{3}}\dots  
\\  = & \ 0
\quad 
\text{because}\  \ell_{2}\ \text{is free in the condition } \ 
2(1+\sum_{t}\ell_{t})  = 1+ \sum_{t} t \ell_{t} . 
\end{align*} 
\kmcomment{
\\\noalign{because \(\ds \ell_{2}\) is free in the condition 
\(\ds 
2(1+\sum_{t}\ell_{t})  = 1+ \sum_{t} t \ell_{t}\) .} 
} 
So, \(\ds \sum_{m>0}(-1)^{m}\myCSW{m}{0,0} = \sum_{m>0} (-1)^{m}
\tbinom{\dim \mathfrakX{1}{1}}{m} = -1\).  When $h=0$, then 0-th chain
space \(\ds \myCSW{0}{0,0} \) is defined and trivially 1-dimensional.
Thus, the Euler number  \(\ds \sum_{m \geq 0}(-1)^{m}\myCSW{m}{0,0} =
0\). 

When $h<0$ then \eqref{eqn:recursive} says that 
\(\ds 
\mynabla{2m+h}{m} = B\cdot \mynabla{2m+h-1}{m-1} \) and we follow the 
same discussion about \(\ds\dim (B\cdot \lambda)\) and get the
conclusion that the Euler number is 0. 

When $h>0$ then \eqref{eqn:recursive} says that 
\(\ds 
\mynabla{2m+h}{m} = B\cdot \mynabla{2m+h-1}{m-1} \sqcup
\T{m} \mynabla{m+h}{m}  
\) and we know the alternating sum is 0 about the first  
term.  Concerning the second term, take an arbitrary element 
\(\ds \lambda = [\ell_{1},\ell_{2},\dots] \in \mynabla{m+h}{m}\) with
the conditions \(\ds \sum_{s} \ell_{s} = m \) and 
\(\ds \sum_{s}s \ell_{s} = m+h \). Then \(\ds \T{m}\cdot \lambda = [0,
\ell{1},\ell_{2},\dots]\) and so 
\begin{align*} 
& \sum_{m} (-1)^m \sum_{\lambda\in\mynabla{m+h}{m}}
\dim( \T{m} \cdot \lambda) 
\\
= &   
\sum_{m} (-1)^m \sum_{\lambda\in\mynabla{m+h}{m}}
\tbinom{\dim \mathfrakX{1}{1}}{\ell{1}}
\tbinom{\dim \mathfrakX{1}{2}}{\ell{2}}\dots 
 =  
\sum_{ 
\sum_{s} \ell_{s} = \sum_{s}s \ell_{s} - h } 
(-1)^{ 
\sum_{s} \ell_{s} 
} 
\tbinom{\dim \mathfrakX{1}{1}}{\ell{1}}
\tbinom{\dim \mathfrakX{1}{2}}{\ell{2}}\dots 
\\  =&  
\sum_{ \ell_{1}} (-1)^{\ell_{1}}  
\tbinom{\dim \mathfrakX{1}{1}}{\ell{1}}
\sum_{ 
\sum_{s} \ell_{s} = \sum_{s}s \ell_{s} - h } 
(-1)^{ 
\sum_{s \ne 1} \ell_{s} 
} 
\tbinom{\dim \mathfrakX{1}{2}}{\ell{2}}
\tbinom{\dim \mathfrakX{1}{3}}{\ell{3}}\dots 
\\ = &  \ 0 \;. 
\end{align*}

\kmqed

\subsubsection{The first weight \(w=1\) case }
Assume $w = 1$. Then using Corollary \ref{kmcor:one} directly, we have
\begin{align*}
\myCSW{m}{1,h} &= \sum_{
 \sum_{s=1}^{m} (h_{s}+1) = h + 2m   
}
\mathfrakX{[m-1,1,0,\dots ]}{(h_1,\dots,h_{m})} 
=  \sum_{ h_{m} }
 \SubCD{1:(m-1)}{ h+1- h_{m} }\mywedge  
 \SubCD{2:1}{ h_{m} } 
 \end{align*}
and 
\(\ds 
 \SubCD{1:(m-1)}{ h+1- h_{m} }\) is just 
\(\ds  \myCSW{m-1}{0, h +1- h_{m} } \).  
Thus, we have next proposition which gives a rule of 
expression of $\ds \myCSW{m}{1,h} $ by lower weight chain spaces 
 $\ds \myCSW{m-1}{0,h'}$.  
\begin{prop} \label{prop:wOne:by:Zero}
The chain complex \(\ds \{ \myCSW{\bullet}{1,h}\} \) is non-trivial if
$\ds h\geq -(1+ \dim \mathfrakX{1}{0}) $, and  
\begin{align}
\myCSW{m}{1,h} =& \sum_{h'}\myCSW{m-1}{0,h-h'+1} 
\mywedge \mathfrakX{2}{h'} \quad\text{for}\quad m \geq  1\;.
\label{align:w1:pos}
\end{align} 
Each degree $m$ of the chain complex is upper bounded by   
\(\ds h+2 + 2\dim 
\mathfrakX{1}{0}+ \dim \mathfrakX{1}{1}\). 

For \(\ds\mR^{n}\) of general $n$, the Euler number of the chain complex 
\(\ds \{ \myCSW{\bullet}{1,h}\} \) is always 0 for each $h$.  
\end{prop}

\textbf{Proof:} 
\eqref{align:w1:pos}
implies 
\(\ds 
\dim \myCSW{m}{1,h} = \sum_{h'}\dim \myCSW{m-1}{0,h-h'+1} \dim
\mathfrakX{2}{h'}\) for $ m \geq  1$.  
\begin{align*}
\sum_{m} (-1)^{m} \dim \myCSW{m}{1,h} = &  
\sum_{m \geq 1} (-1)^{m} 
\sum_{h'}\dim \myCSW{m-1}{0,h-h'+1} \dim
\mathfrakX{2}{h'} \\ = & 
- 
\sum_{h'} \dim \mathfrakX{2}{h'} 
\sum_{m \geq 1} (-1)^{m-1} \dim \myCSW{m-1}{0,h-h'+1} 
\\ = & \quad 
 0 \quad \text{using  
 Lemma \ref{thm:Euler:wZero}.}  
\end{align*}
\kmqed

\kmcomment{
\begin{exam} We apply the Proposition above in the case of $n=2$ and
$w=1$ to know the dimension of each chain space and also Euler number
for each $h$. We only show the tables below in lower $h$.  Euler number
is zero for the cases below. 

$w=1, h=-3$
\begin{center}
\begin{tabular}{c | *{7}{c} }
m      & 1 & 2 & 3  & 4  & 5  & 6  & 7 \\\hline
$\dim$ & 0 & 0 & 1 & 4 & 6 & 4 & 1  
\end{tabular}
\end{center} 
\kmcomment{
$w=1, h=-2$
\begin{center}
\begin{tabular}{c | *{8}{c} }
m      & 1 & 2 & 3  & 4  & 5  & 6  & 7  & 8 \\\hline
$\dim$ & 0 & 2 & 10 & 26 & 44 & 46 & 26 & 6 
\end{tabular}
\end{center}
$w=1, h=-1$
\begin{center}
\begin{tabular}{c | *{9}{c} }
m      & 1 & 2 & 3  & 4  & 5  & 6  & 7  & 8 & 9 \\\hline
$\dim$ & 1 & 8 & 37 & 108 & 202 & 244 & 185 & 80 & 15 
\end{tabular}
\end{center}
}
$w=1, h=0$
\begin{center}
\begin{tabular}{c | *{10}{c} }
m      & 1 & 2 & 3  & 4  & 5  & 6  & 7  & 8 & 9  & 10\\\hline
$\dim$ & 2 & 20 & 104 & 330 & 688 & 964 & 888 & 506 & 158 & 20
\end{tabular}
\end{center}
\kmcomment{
$w=1, h=1$
\begin{center}
\begin{tabular}{c | *{11}{c} }
m      & 1 & 2 & 3  & 4  & 5  & 6  & 7  & 8 & 9  & 10 & 11\\\hline
$\dim$ & 3 & 40 & 238 & 848 & 1976 & 3112 & 3321 & 2332 & 999 & 220 & 15 
\end{tabular}
\end{center}\small
}
$w=1, h=2$
\begin{center}
\begin{tabular}{c | *{12}{c} }
m      & 1 & 2 & 3  & 4  & 5  & 6  & 7  & 8 & 9  & 10 & 11 & 12\\\hline
$\dim$ & 4 & 70 & 484 & 1924 & 4968 & 8734 & 10570 & 8670 & 4608 & 1444
& 214 & 6 
\end{tabular}
\end{center}
\end{exam}

Combining  Lemma \ref{thm:Euler:wZero} and 
Proposition \ref{prop:wOne:by:Zero}, 
we have next result of  
the Euler number for $w=1$ case: 
\begin{prop} \label{prop:w1}
For \(\s\mR^{n}\) of general $n$, the Euler number of chain complex 
\(\ds \{ \myCSW{\bullet}{1,h}\} \) is always 0 for each $h$.  
\end{prop}
\textbf{Proof:} 
Proposition \ref{prop:wOne:by:Zero} says  
\(\ds 
\dim \myCSW{m}{1,h} = \sum_{h'}\dim \myCSW{m-1}{0,h-h'+1} \dim
\mathfrakX{2}{h'}\) for $ m > 1$ and 
\(\ds \dim \myCSW{1}{1,h} = \dim \mathfrakX{2}{h+1}\). 
\begin{align*}
\sum_{m} (-1)^{m} \dim \myCSW{m}{1,h} = &  
- \mathfrakX{2}{h+1} + \sum_{m>1} (-1)^{m} 
\sum_{h'}\dim \myCSW{m-1}{0,h-h'+1} \dim
\mathfrakX{2}{h'} \\ = & 
- \mathfrakX{2}{h+1} - 
\sum_{h'} \dim \mathfrakX{2}{h'} 
\sum_{m>1} (-1)^{m-1} 
\dim \myCSW{m-1}{0,h-h'+1} 
\\ = & 
- \mathfrakX{2}{h+1} + 
 \mathfrakX{2}{h+1} = 0 \quad \text{using  
 Lemma \ref{thm:Euler:wZero}}\;.  
\end{align*}
\kmqed
}

\subsubsection{The first weight \(w=2\) case }
Assume $w = 2$. Again, using Corollary \ref{kmcor:one}, we have
\begin{align*}
\myCSW{m}{2,h} &= 
\sum_{ 
 \sum_{s=1}^{m} (h_{s}+1) = h + 2m}   
\mathfrakX{[m-1,0,1,0,\dots ]}{(h_1,\dots,h_{m})}
+ \sum _{
 \sum_{s=1}^{m} (h_{s}+1) = h + 2m}   
\mathfrakX{[m-2,2,0,\dots ]}{(h_1,\dots,h_{m})}
\\ 
& = \sum 
 \SubCD{1:(m-1)}{ h +1 - h_{m}} \mywedge \mathfrakX{3}{ h_{m} } 
 + \SubCD{1:(m-2)}{ h -h' }\mywedge\SubCD{2:2}{h'} \;. 
\end{align*} 
Thus, we have the next proposition which gives a rule of 
expression of $\ds \myCSW{m}{2,h} $ by lower weight chain spaces.  
\begin{prop}
\label{prop:wTwo:by:Zero} 
The chain complex \(\ds \{ \myCSW{\bullet}{2,h}\} \) is non-trivial if
$\ds h\geq -(2+ \dim\mathfrakX{1}{0}) $, and  
\begin{align}
\myCSW{m}{2,h} =& \sum_{h'}\myCSW{m-1}{0,h+1-h'} 
\mywedge \mathfrakX{3}{h'} 
\sqcup \sum_{ a\leq b} \myCSW{m-2}{0,h+2-a-b}
\mywedge \mathfrakX{2}{a} \mywedge \mathfrakX{2}{b} 
\quad\text{for}\quad m \geq  2\;,  \label{gen:two:h}
\\
\kmcomment{
\myCSW{2}{2,h} =&
\sum_{ h_{1}+h_{2} = 2+ h }
\mathfrakX{1}{h_{1}}\mywedge \mathfrakX{3}{h_{2}}
\sqcup  \sum_{ \substack{ h_{1}+h_{2} = 2+ h \\ h_{1}\leq h_{2} } }
\mathfrakX{2}{h_{1}}\mywedge \mathfrakX{2}{h_{2}} 
\label{two:two:h} \;, \\ 
}
\noalign{and} 
\myCSW{1}{2,h} =& \mathfrakX{3}{h+1}\;.  \label{one:two:h} 
\end{align} 
The range of degree $m$ of the chain complex has an upper bound 
\(\ds h+4 + 2\dim \mathfrakX{1}{0}+ \dim \mathfrakX{1}{1}\).  
\end{prop}


We can apply Lemma \ref{thm:Euler:wZero} for the chain complex $w=2$, 
we have 

\begin{prop} \label{prop:w2}
For general $n$, the Euler number of chain complex 
\(\ds \{ \myCSW{\bullet}{2,h}\} \) is always 0 for each $h$.  
\end{prop}
\textbf{Proof:}
We only alternating sum up the  dimension of chain spaces. 
\kmcomment{
Since first term of 
\eqref{gen:two:h} is the same with 
\(\ds
\sum_{ h' }
\myCSW{1}{0,h+1-h'} 
\mywedge \mathfrakX{3}{h'}
\),}
We first sum up the terms which involve \(\ds \mathfrakX{3}{\bullet}\) as follows:
\begin{align*}
A =&  (-1)^1 \dim \mathfrakX{3}{h+1}  
+ \sum_{m \geq 2} (-1)^{m} 
\sum_{h'} \dim \myCSW{m-1}{0,h+1-h'} 
\dim \mathfrakX{3}{h'} 
\\ =& 
- \sum_{ h' }\dim \mathfrakX{3}{h'} 
 \sum_{m \geq 0} (-1)^{m} 
\dim \myCSW{m}{0,h+1-h'} 
=  0 \quad \text{using Lemma \ref{thm:Euler:wZero}. 
}
\end{align*}
The rest is 
\begin{align*}
B =& 
\sum_{m \geq 2} (-1)^{m} 
\sum_{ a\leq b} \dim \myCSW{m-2}{0,h+2-a-b}
\dim (
 \mathfrakX{2}{a} \mywedge \mathfrakX{2}{b} )  \\
 =& 
\sum_{ a\leq b} 
\dim ( \mathfrakX{2}{a} \mywedge \mathfrakX{2}{b} )  
\sum_{m\geq 2} (-1)^{m-2} 
\dim \myCSW{m-2}{0,h+2-a-b}
=  0 \quad  
\text{using again Lemma \ref{thm:Euler:wZero}.} 
\end{align*}
\kmqed

\subsubsection{
General first weight case 
}
Inspired by Propositions \ref{prop:wOne:by:Zero}
and \ref{prop:w2}, we have the next general result including those results.   
\begin{thm}\label{thm:triv:general}
For general $n$, the Euler number of chain complex 
\(\ds \{ \myCSW{\bullet}{w,h}\} \) is 0 for 
each $w$ and $h$.  
\end{thm}
\textbf{Proof:} 
We have already seen that it is true for $w=0,1,2$. 
So we may assume $w>2$ and $m>0$.    
We use the notation \eqref{eqn:new:notation}. From Corollary
\ref{kmcor:one}, we have the chain space is written by 
\[\ds
 \myCSW{m}{w,h} = \mathop{\oplus}_{
\substack{
\sum_{i=1}^{n} k_{i} = m \\ 
\sum_{i=1}^{n} i k_{i} = w+m\\
\sum_{i=1}^{n}  u_{i} = h 
}
}
\SubCD{ 1:k_{1}} {u_{1}} \mywedge\cdots\mywedge 
\SubCD{ n:k_{n}} {u_{n}}\;. \] 
 The first component 
\(\ds \SubCD{ 1:k_{1}} {u_{1}}\) is equal to the chain space \(\ds 
 \myCSW{k_{1}}{0,u_{1}}\) with the first weight 0. Thus 
\begin{align*}
 & \sum_{m>0} (-1)^m \dim \myCSW{m}{w,h}\\
=& 
\sum_{ \sum_{i=1}^{n} (i-1) k_{i} = w }
(-1)^{ \sum_{s=1}^{n} k_{s} } \sum_{ u_{j}} \dim 
 \myCSW{k_{1}}{0,u_{1}} \dim\big( 
\SubCD{ 2:k_{2}} {u_{2}} \mywedge\cdots\mywedge 
\SubCD{ n:k_{n}} {u_{n}}\big) \\
=& 
\sum_{ \sum_{i=1}^{n} (i-1) k_{i} = w }
(-1)^{ \sum_{s=2}^{n} k_{s} } \sum_{ u_{j}}\sum_{k_{1}}
(-1)^{ k_{1} } \dim 
 \myCSW{k_{1}}{0,u_{1}} \dim\big( 
\SubCD{ 2:k_{2}} {u_{2}} \mywedge\cdots\mywedge 
\SubCD{ n:k_{n}} {u_{n}}\big)
\\\noalign{we used that the condition \(\ds   
 \sum_{i=1}^{n} (i-1) k_{i} = w \) does not involve \(\ds k_{1}\), and  
now we use Lemma \ref{thm:Euler:wZero}}
=& \sum_{ \sum_{i=1}^{n} (i-1) k_{i} = w }
(-1)^{ \sum_{s=2}^{n} k_{s} } \sum_{ u_{j}} 0 
=  0 \;. 
\end{align*}
\kmqed

\subsubsection{Extended case}
 \(\ds \mathfrakX{0}{} (M) \oplus \mathfrakX{1}{} (M) \oplus
\cdots \oplus \mathfrakX{{n}}{} (M) \) is also 
a pre Lie superalgebra including  
\(\ds \mathfrakX{1}{} (M) \oplus
\cdots \oplus \mathfrakX{{n}}{} (M) \) which we dealt so far.

In this subsection again taking \(\ds M = \mR^{n}\) we consider the chain spaces defined by 
\[\myOSW{m+1}{w,h} = \sum_{ } 
\mathfrakX{0}{h_{0}} (\mR^{n}) \mywedge
\mathfrakX{i_{1}}{h_{1}} (\mR^{n}) \mywedge
\cdots \mywedge 
\mathfrakX{i_{m}}{h_{m}} (\mR^{n}) \] 
where \(\ds (0-1) + \sum_{s=1}^{m} (i_{s}-1) = w\) and  
\(\ds \sum_{s=0}^{m} (h_{s}-1) = h\). 
We see easily the next proposition. 
\begin{prop}
\(\ds \pdel 
(\myOSW{m+1}{w,h}) \subset 
\myOSW{m}{w,h} \) and have another double-weighted homology groups. 
The Euler number of the chain complex  
\( \ds \{ \myOSW{\bullet}{w,h} \} \) is 0 . 
\end{prop}

\subsection{Homology groups with non-trivial representation} 
In this subsection, we first study natural representation of 
 \(\ds \frakg = \mathfrakX{1}{} (M) \oplus
\cdots \oplus \mathfrakX{{n}}{} (M) \) for general manifold $M$.  
Since the Schouten bracket of  \(\ds \frakg \) and 
\(\ds \mathfrakX{0}{} (M) = C^{\infty}(M) \) lies in 
 \(\ds \mathfrakX{0}{} (M) \oplus
\cdots \oplus 
\mathfrakX{{n-1}}{} (M) \), we regard \(\ds \frakg \) acts on  
\(\ds \mathfrakX{0}{} (M) = C^{\infty}(M) \) by 
\[ U \cdot f = \Sbt{U}{f} \mod \frakg \;. \]
Actually, if \(\ds U \in \mathfrakX{1}{} (M) \) then \(\ds U\cdot f =
\Sbt{U}{f} = \langle U, df\rangle\) and 
if \(\ds U \in \mathfrakX{i}{} (M) \) then \(\ds U\cdot f = 0\) for $i>1$. 
So we have a representation space \(\ds V= 
\mathfrakX{0}{} (M) = C^{\infty}(M) \) of \(\ds \frakg\) and the action.   

Now we consider \(\ds M = \mR^{n}\) and we may study relative 
homology groups of the chain spaces 
\(\ds \Delta^{m} \frakg \otimes V \) with 
the boundary operator \(\ds \pdel_{V}\) as
introduced in the section \ref{sec:prelim}. 

Now we introduce double-weighted chain spaces using the specialty 
of our base space \(\ds M = \mR^{n}\). 
The chain spaces are given by 
\kmcomment{
\myOSW{0}{w,h} &=
\mathfrakX{0}{1+h} (\mR^{n}) \quad (\forall w)\;, \label{non:triv:zero}
\\
}
\begin{align} 
\myOSW{m}{w,h} 
&= \sum_{\substack{
 \sum_{s=1}^{m} (i_{s}-1) = w\\ 
 \sum_{s=0}^{m} (h_{s}-1) = h 
}} 
\mathfrakX{i_{1}}{h_{1}} (\mR^{n}) \mywedge
\cdots \mywedge 
\mathfrakX{i_{m}}{h_{m}} (\mR^{n}) 
\otimes 
\mathfrakX{0}{h_{0}} (\mR^{n})  \quad (m \geq 0)\;  \label{non:triv:one}\\
&= \sum_{h_{0}} 
\myCSW{m}{w, h+1-h_{0} } 
\otimes 
\mathfrakX{0}{h_{0}} (\mR^{n})  \quad (m \geq 0)\;,  \label{non:triv:two}
\end{align} 
where \(\ds 
\myCSW{\bullet}{w, h'} \) are the chain spaces in the trivial module. 
We easily see the next proposition. 
\begin{prop}
The double-weight is invariant by 
\(\ds \pdel_{V}\), i.e., \(\ds \pdel_{V} \left(
\myOSW{m}{w,h}\right) \subset 
\myOSW{m-1}{w,h}\). Thus, we have the double-weighted homology groups
\(\ds \myHomW{m}{w,h}(\frakg, V) \)  with $\frakg$-module \(V\).   
\end{prop}
\textbf{Proof:} 
We know that \(\ds \Sbt{ \mathfrakX{i}{h}} { \mathfrakX{i'}{h'}} \subset
\mathfrakX{i+i'-1}{h+h'-1} \) and in particular, 
\(\ds \Sbt{ \mathfrakX{i}{h}} { \mathfrakX{0}{h'}} \subset
\mathfrakX{i-1}{h+h'-1} \).
\begin{align*}
\pdel_{V} ( 
\mathfrakX{i_{1}}{h_{1}}  \mywedge
\cdots \mywedge 
\mathfrakX{i_{m}}{h_{m}}  
\otimes 
\mathfrakX{0}{h_{0}}  )  = & 
\pdel (
\mathfrakX{i_{1}}{h_{1}}  \mywedge
\cdots \mywedge \mathfrakX{i_{m}}{h_{m}}  ) \otimes 
\mathfrakX{0}{h_{0}}   \pm  
\sum_{p} \mathfrakX{i_{1}}{h_{1}} \Delta \dots
\widehat{\mathfrakX{i_{p}}{h_{p}}} \dots {\mathfrakX{i_{m}}{h_{m}}}
\otimes \Sbt{ \mathfrakX{i_{p}}{h_{p}}}  {\mathfrakX{0}{h_{0}}}
\end{align*} 
Thus, we directly see that 
the double-weight of the first part 
\(\ds  \Sbt{ \mathfrakX{i_{p}}{h_{p}}} { \mathfrakX{i_{q}}{h_{q}}} \Delta
\dots \widehat{\mathfrakX{i_{p}}{h_{p}}} \dots \widehat 
{\mathfrakX{i_{q}}{h_{q}}}  \dots \otimes 
 \mathfrakX{0}{h_{0}} \) does not change. 
About the second part,  
\(\ds \mathfrakX{i_{1}}{h_{1}} \Delta \dots \widehat{\mathfrakX{i_{p}}{h_{p}}} \dots 
{\mathfrakX{i_{m}}{h_{m}}}  \otimes \Sbt{ 
\mathfrakX{i_{p}}{h_{p}}}  
{\mathfrakX{0}{h_{0}}} \) is  0  if \(\ds i_{p} \ne 1\). 
When  \(\ds i_{p} = 1\), the first weight is \( \ds 
\sum_{s\ne p} (i_{s}-1)  = 
\sum_{s=1}^{m} (i_{s}-1) = w \) and    
the second weight is \( \ds 
\sum_{s\ne p} (h_{s}-1) + (h_{p}+h_{0} -1-1) = 
\sum_{s=0}^{m} (h_{s}-1) = h \).   
\kmqed

\medskip
Due to Lemma \ref{thm:Euler:wZero} and 
Theorem \ref{thm:triv:general}, we have  
next result about Euler number of the chain complex 
\(\ds (\myOSW{\bullet}{w,h},  \pdel_{V}) \). 
\begin{thm}\label{thm:module:gen}
The Euler number of  
\(\ds (\myOSW{\bullet}{w,h},  \pdel_{V}) \) is 0 for each $w$ and $h$.   
\kmcomment{
If $w>0$, then 
\(\ds \sum_{m>0} (-1)^{m} \dim (\myOSW{m}{w,h})\) is 0, and  
the Euler number of  
\(\ds (\myOSW{\bullet}{w,h},  \pdel_{V}) \) is \(\ds
\tbinom{n+h}{n-1}\). }
\end{thm}
\textbf{Proof:}
From \eqref{non:triv:two}, 
\begin{align*} 
\sum_{m\geq 0} (-1)^{m} \dim \myOSW{m}{w,h} = &  
\sum_{m\geq 0} (-1)^{m} \sum_{h_{0}} \dim \myCSW{m}{w,h+1-h_{0}} \dim
\mathfrakX{0}{h_{0}} \\ =  & 
\sum_{h_{0}} \dim \mathfrakX{0}{h_{0}}    \sum_{m\geq 0} (-1)^{m} 
\dim \myCSW{m}{w,h+1-h_{0}} \\\noalign{using Theorem
\ref{thm:triv:general} }
 = &  
\sum_{h_{0}} \dim \mathfrakX{0}{h_{0}} \times  0 
= 0\;.     
\end{align*}
\kmqed

\subsection{Example of pre Lie superalgebra related to a Lie
superalgebra}
In Example \ref{exam:super:elem}, we saw toy models of Lie superalgebra
and pre Lie superalgebra. We study the chain complexes of those. 

\(\ds \frakg = \mathfrak{gl}(2) = 
\frakg_{0} \oplus \frakg_{1}  \oplus \frakg_{2}\). Take a basis 
\( \futoji{u}{1} \in \frakg_{0}\), 
\( \futoji{u}{2},\futoji{u}{3} \in \frakg_{1}\), 
\( \futoji{u}{4} \in \frakg_{2}\) with the following bracket relation: 
\begin{center}
\begin{tabular}{c|*{4}{c}} 
 & \(\futoji{u}{1}\) & \(\futoji{u}{2}\) & \(\futoji{u}{3}\) & 
 \(\futoji{u}{4}\) \\\hline
 \(\futoji{u}{1}\) & 0 & $ 2\futoji{u}{2}$ & $-2 \futoji{u}{3}$ & 0 \\\hline
 \(\futoji{u}{2}\) & $-2 \futoji{u}{2}$ & 0 & $\futoji{u}{4}$ & 0 \\\hline
 \(\futoji{u}{3}\) & $2\futoji{u}{3}$ & $\futoji{u}{4}$ & 0 & 0 \\\hline
 \(\futoji{u}{4}\) & 0 & 0 & 0 & 0 \\\hline
\end{tabular}
\hfil 
\begin{tabular} { *{4}{c}}
 m & w-1 & w & w+1 \\\hline
\(\ds\myCSW{m}{w}\) & \(\ds \Delta^{w-2} \frakg_{1} \Delta\frakg_{2}\) &
\(\ds \Delta^{w} \frakg_{1} \oplus \frakg_{0}
\Delta^{w-2}\frakg_{1}\Delta\frakg_{2}\) & \(\ds \frakg_{0} \Delta ^{w}
\frakg_{1}\) \\
$\dim$ & w-1 & 2w& w+1 \\
$\dim \pdel$ & w-1 & w+1 & 0 \\\hline
Betti & 0 &  0 & 0 \end{tabular} 
\end{center}
Given a weight $w$, the chain spaces are given by \(\ds\myCSW{m}{w}\)
for $m=w-1, w, w+1$ and we get dimension, rank and Betti numbers as
right above.  

\medskip

Suppose 
\(\ds i' = i + I p \) and 
\(\ds j' = j + J p \), i.e.,  \( i \equiv i'\) and  \( j \equiv j'\)
\(\mod  p\).  Then 
\[ i'+ j' = i+j+ (I+J) p \equiv i+j \mod p \] and   
\[ i'  j' = (i+I p) (j+ J p) = i j + i J p + j I p + I J p^2 \equiv i j
\mod p
\]  but 
\(\ds i' j' - i j  =  i J p + j I p + I J p^2 \equiv 0 
\mod 2 
\) if $p$ is even.  So Lie super algebras should be divided into even
number of subspaces.  Thus, our example 
\(\ds \frakg = \mathfrak{gl}(2) = 
\frakg_{0} \oplus \frakg_{1}  \oplus \frakg_{2}\) is not Lie
superalgebra. 

The Lie superalgebra 
in Example \ref{exam:super:elem} is sometimes denoted as 
\(\ds \mathfrak{gl}(1|1) = \frakg_{[0]} \oplus \frakg_{[1]} \). 

\(\ds \frakg_{[0]} \) is spanned by 
  \(\futoji{u}{1}\) and  \(\futoji{u}{2}\), and  
\(\ds  \frakg_{[1]} \) is spanned by  
  \(\futoji{u}{3}\) and  \(\futoji{u}{4}\).  Those basis satisfy the
  next bracket relations:   
\begin{center}
\begin{tabular}{c|*{4}{c}} 
 & \(\futoji{u}{1}\) & \(\futoji{u}{2}\) & \(\futoji{u}{3}\) & 
 \(\futoji{u}{4}\) \\\hline
 \(\futoji{u}{1}\) & 0 & 0 & $\futoji{u}{3}$ & $-\futoji{u}{4}$\\\hline
 \(\futoji{u}{2}\) & 0 & 0 & $-\futoji{u}{3}$ & $\futoji{u}{4}$\\\hline
 \(\futoji{u}{3}\) & $-\futoji{u}{3}$ & $\futoji{u}{3}$ & 0 & 
  $\futoji{u}{1} + \futoji{u}{2}$  \\\hline
 \(\futoji{u}{4}\) & $\futoji{u}{4}$ & $-\futoji{u}{4}$ &  
  $\futoji{u}{1} + \futoji{u}{2}$ & 0  \\\hline 
\end{tabular}
\end{center}

Depending on the weight to be even or odd, we have 
\begin{align*}
 \myCSW{2m}{[0]} &= 
 \Delta^{2m} \frakg_{[1]} \oplus 
\Delta^{2} \frakg_{[0]} \Delta^{2m-2} \frakg_{[1]}\;, &
\myCSW{2m+1}{[0]} &= 
\frakg_{[0]} \Delta^{2m} \frakg_{[1]}\;,\\ 
 \myCSW{2m}{[1]} &= 
\frakg_{[0]} \Delta^{2m-1} \frakg_{[1]}\;, &
\myCSW{2m+1}{[1]} &= 
 \Delta^{2m+1} \frakg_{[1]} \oplus 
\Delta^{2} \frakg_{[0]} \Delta^{2m-1} \frakg_{[1]}
\;.  
\end{align*} 

Denote \(\ds \Delta^{a} \futoji{u}{3} 
\Delta^{b} \futoji{u}{4} \) by \(F(a,b)\). Then 
\begin{align*} 
 \pdel (F(a,b)) & = a b (  \futoji{u}{1}+ \futoji{u}{2})\Delta F(a-1,b-1) \\ 
 \pdel (\futoji{u}{1}\Delta F(a,b)) & = 
- a b   \futoji{u}{1} \Delta  \futoji{u}{2}\Delta F(a-1,b-1) 
+ (a- b )  F(a,b) 
\\ 
\pdel (\futoji{u}{2}\Delta F(a,b)) & = 
 a b   \futoji{u}{1} \Delta  \futoji{u}{2}\Delta F(a-1,b-1) 
- (a- b )  F(a,b) 
\\ 
 \pdel (\futoji{u}{1}\Delta 
\futoji{u}{2}\Delta 
F(a,b)) & = (a-b) \futoji{u}{1}\Delta \futoji{u}{2}\Delta F(a,b)
\end{align*}

\begin{center}
\begin{tabular}{c | *{5}{c} }
$\bullet $ & \dots  & $2m-1$ & $2m$ & $2m+1$  \\\hline
$\dim \myCSW{\bullet}{[0]}$ & \dots & $2(2m-1)$ & $2(2m)$ & $2(2m+1)$ \\
$\dim\pdel$ & $2m-1$ & $2m-1$ & $2m+1$ & $2m+1$ \\\hline
Betti & \dots  &  0 &  0 &  0 
\end{tabular}
\\[3mm] 
\begin{tabular}{c | *{5}{c} }
$\bullet $ & $2m-2$ & $2m-1$ & $2m$ & $2m+1$    \\\hline
$\dim \myCSW{\bullet}{[1]}$ & $2(2m-2)$ & $2(2m-1)$ & $2(2m)$ & $2(2m+1)$ \\
$\dim\pdel$ & $2(m-1)$ & $2m$ & $2m$ & $2(m+1)$ \\\hline
Betti & \dots  &  0 &  0 &  0 
\end{tabular}
\end{center}

\kmcomment{
Since 
\begin{align*}
 \myCSW{m}{[0]} =& \sum_{ \substack{ 0\leq m-2k\leq 2 } }
\Delta^{m-2k} \frakg_{[0]} \Delta^{2k} \frakg_{[1]}\;, \quad 
\myCSW{m}{[1]} = \sum_{ \substack{ 0\leq m-1-2k\leq 2} }
\Delta^{m-1-2k} \frakg_{[0]} \Delta^{2k+1} \frakg_{[1]}\;, 
\end{align*} 
\(\ds \dim 
 \myCSW{m}{[0]} = \sum_{ 
k} 
 \tbinom{\dim \frakg_{[0]} }{ m-2k}
 \tbinom{\dim \frakg_{[1]}-1 + 2k} {\dim \frakg_{[1]}-1}
\) and we expect 
\begin{align*}
\sum_{m} (-1)^{m} 
 \dim \myCSW{m}{[0]} 
= &  \sum_{m} (-1)^{m}  \sum_{
k
 } 
 \tbinom{\dim \frakg_{[0]} }{ m-2k}
 \tbinom{\dim \frakg_{[1]}-1 + 2k} {\dim \frakg_{[1]}-1}
 \\ 
\mathop{=}^{?} & \sum_{k} 
 \tbinom{\dim \frakg_{[1]}-1 + 2k} {\dim \frakg_{[1]}-1}
 \sum_{m} (-1) ^{m} 
 \tbinom{\dim \frakg_{[0]} }{ m-2k} = 0 \;? \end{align*}
Also we have the same  expectation to  
\(\ds \{\myCSW{m}{[1]}\}_{m} \).  

\begin{align*}
 \myCSW{m}{[0]} =& \sum_{
\substack{ 0\leq m-2k\leq \dim \frakg_{[0]} } 
 }
\Delta^{m-2k} \frakg_{[0]} \Delta^{2k} \frakg_{[1]}\;, \quad 
\myCSW{m}{[1]} = \sum_{
\substack{ 0\leq m-1-2k\leq \dim \frakg_{[0]} } 
}
\Delta^{m-1-2k} \frakg_{[0]} \Delta^{2k+1} \frakg_{[1]}\;. \end{align*} 
Since \(\ds \dim 
 \myCSW{m}{[0]} = \sum_{
k
 } 
 \tbinom{\dim \frakg_{[0]} }{ m-2k}
 \tbinom{\dim \frakg_{[1]}-1 + 2k} {\dim \frakg_{[1]}-1}
\), we expect 
\begin{align*}
\sum_{m} (-1)^{m} 
 \dim \myCSW{m}{[0]} 
= &  \sum_{m} (-1)^{m}  \sum_{
k
 } 
 \tbinom{\dim \frakg_{[0]} }{ m-2k}
 \tbinom{\dim \frakg_{[1]}-1 + 2k} {\dim \frakg_{[1]}-1}
 \\ 
\mathop{=}^{?} & \sum_{k} 
 \tbinom{\dim \frakg_{[1]}-1 + 2k} {\dim \frakg_{[1]}-1}
 \sum_{m} (-1) ^{m} 
 \tbinom{\dim \frakg_{[0]} }{ m-2k} = 0 \;? \end{align*}
Also we have the same  expectation to  
\(\ds \{\myCSW{m}{[1]}\}_{m} \).  
}

\section{Betti numbers of homology groups of concrete pre Lie
superalgebras}\label{sec:Betti} 
So far, we studied the chain spaces \(\myCSW{m}{w,h}\) for fixed space
dimension $n$ and double weight \((w,h)\).  
As stated in Remark \ref{remark:Poisson}, we may find all Poisson
structures in the second homology group of pre Lie superalgebra of
tangent bundle of $M$ with the Schouten bracket. Thus, it is interesting
to study the second and/or the third  homology group.   
But, it seems hard to attack to
general manifold $M$.  So, again we deal with 
the pre Lie superalgebra of 
homogeneous polynomial coefficients multi vector fields on \(\ds
\mR^{n}\). In this section, we study not only the second Betti number
but also Betti numbers of general degree.

In general pre Lie superalgebra theory, recursive formulae of the boundary operator is given
as below in two ways: one is given by using right action and the other is given
by left action. 
\begin{align} 
\pdel (A_{1} \mywedge \cdots \mywedge A_{m+1}) &= 
\pdel (A_{1} \mywedge \cdots \mywedge A_{m}) \mywedge A_{m+1} 
+ (-1)^{m+1} 
(A_{1} \mywedge \cdots \mywedge A_{m})^{ A_{m+1} }
\label{pdel:recurs:right}
\\\noalign{where} \label{rigtht:action}
(A_{1} \mywedge \cdots \mywedge A_{m})^{ A_{m+1} } &= 
 (A_{1} \mywedge \cdots \mywedge A_{m-1}) \mywedge \Sbt{A_{m}}{ A_{m+1} }  
\\ & \qquad  
+ 
(-1)^{a_{m}a_{m+1}}
(A_{1} \mywedge \cdots \mywedge A_{m-1})^{ A_{m+1} } \mywedge A_{m} 
\notag
\\
&= \sum_{i=1}^{m} (-1)^{a_{m+1} \sum_{s=i+1}^{m} a_{s} } A_{1} \mywedge \cdots 
\mywedge \Sbt{ A_{i} }{A_{m+1}} \mywedge \cdots \mywedge A_{m} 
\\ 
\label{pdel:recurs:left}
\pdel (A_{0} \mywedge A_{1} \mywedge 
\cdots \mywedge A_{m}) &= 
- A_{0} \mywedge \pdel (A_{1} \mywedge \cdots \mywedge A_{m})
+ A_{0} \cdot  (A_{1} \mywedge \cdots \mywedge A_{m})
\\\noalign{where} \label{left:action}
A_{0} \cdot  (A_{1} \mywedge \cdots \mywedge A_{m}) &=
\Sbt{A_{0}}{A_{1}} \mywedge  
(A_{2} \mywedge \cdots \mywedge A_{m})
+ (-1)^{a_{0}a_{1}} A_{1} \mywedge \left(
A_{0} \cdot  (A_{2} \mywedge \cdots \mywedge A_{m})
\right)  \\
& = \sum_{i=1}^{m} (-1)^{ a_{0}\sum_{s<i}a_{s}} 
A_{1} \mywedge \cdots \mywedge  \Sbt{A_{0}}{A_{i}} \mywedge \cdots   
\mywedge A_{m}) 
\end{align}
for each homogeneous elements \( A_{i}\in\frakg_{a_{i}}\).  
In lower degree, the boundary operator is given
as bellows:  
\begin{align} 
\pdel ( A \mywedge B ) &= \Sbt{A}{B} \\
\pdel ( A \mywedge B \mywedge C ) &= -  {A} \mywedge {
 \Sbt{B}{C} } + 
 \Sbt{A}{B} \mywedge C + (-1)^{a b} B \mywedge  
 \Sbt{A}{C}  
\end{align}
for each homogeneous elements 
\( A\in\frakg_{a}\), \( B\in\frakg_{b}\), \( C\in\frakg_{c}\).

\kmcomment{
In the case of tangent bundle with Schouten bracket, 
we  have the following expression: 
Let \(X,Y \in \tbdl{M}\), \(\pi,\pi' \in \Lambda^{2}\tbdl{M}\), 
\( U \in \Lambda^{3}\tbdl{M}\).  Then we have
\begin{align}
\pdel ( X \mywedge Y \mywedge U  ) =& 
- X \mywedge \Sbt{Y}{U} + \Sbt{X}{Y} \mywedge U 
+ Y \mywedge \Sbt{X}{U} \\ 
\pdel ( X \mywedge \pi \mywedge \pi'  ) =& 
- X \mywedge \Sbt{\pi}{\pi'} + 
\Sbt{X}{\pi} \mywedge \pi' 
+ \Sbt{X}{\pi'}  \mywedge \pi
\end{align}
}
If we will handle Poisson structures on \(\ds\mR^{n}\) by homology
theory of pre Lie superalgebra, then  
Remark \ref{remark:Poisson} 
says 
we will deal with \(\ds\{ \myCSW{\bullet}{w=2,h}\}\), where   
\begin{align*}\myCSW{1}{2,h} & = \mathfrakX{3}{h+1} \;, \\
\myCSW{2}{2,h} & = 
\sum_{a+b=h+2} \mathfrakX{1}{a}\mywedge \mathfrakX{3}{b}
+ \sum_{a+b=h+2} \mathfrakX{2}{a}\mywedge \mathfrakX{2}{b}\;, 
\\
\myCSW{3}{2,h} & = \sum_{c+a+b=h+2+1}
\mathfrakX{1}{c}\mywedge \mathfrakX{1}{a}\mywedge \mathfrakX{3}{b}
+ \sum_{c+a+b=h+2+1} \mathfrakX{1}{c}\mywedge
\mathfrakX{2}{a}\mywedge \mathfrakX{2}{b} \;,
\\ & \hspace{-5mm} \vdots 
\end{align*}
\begin{remark}
Since 
\(\ds
\myCSW{m}{w,h}  = \sum_{\substack{ \sum_{i=1}^{m} a_{i} = w+m,  
\\
\sum_{i=1}^{m} b_{i} = h+m} }
\mathfrakX{a_{1}}{b_{1}}\mywedge 
\mathfrakX{a_{2}}{b_{2}}\mywedge \cdots \mywedge 
\mathfrakX{a_{m}}{b_{m}}\) in general, if 
\(\ds
\myCSW{m}{w,h}  \ne (0)\) then \(\ds a_{i} \leq n\) 
 for each $i$, and so 
 \( \sum_{i=1}^{m}a_{i} \leq mn\). Thus,  
 \(\ds w \leq m(n-1)\). Namely, $w$ is bounded from above by the dimension
 $n$ and the degree $m$ of the chain space.  
\end{remark}

Since \(\ds M = \mR^{n}\), we have a special vector field \( \XO =
\sum_{i=1}^n x_{i} \frac{\pdel}{\pdel x_{i}} \in \mathfrakX{1}{1} \)
(called Euler vector field).  It is known that 
if $f$ is $h$-homogeneous polynomial, then \(\Sbt{\XO}{f} = h f\). 
If \(\ds D \in \mathfrakX{p}{0}\), then \(\Sbt{\XO}{D} = - p D\). 
We have the next lemma in general: 
\begin{lemma} \label{lemma:dilation}
For each $\ds U\in \mathfrakX{p}{h}$, 
\begin{align} 
& \XO \cdot U = 
\Sbt{\XO}{U} = (-p+h)U \\
\noalign{holds.  In fact, the action of \(\XO\) is divided into two parts:}
& \sum_{k=1}^{n} x_{k} \Sbt{ \bibun{k} }{U} = h U \;, 
\\\noalign{and}
& \sum_{k=1}^{n} \bibun{k} \wedge \Sbt{ x_{k} }{U} = -p U\;. 
\\\noalign{Thus,}
& \XO\cdot W = (-w+h) W\quad (\forall W \in \myCSW{m,w,h})\;. 
\label{lemma:chain:map}
\end{align} 
\end{lemma}
\textbf{Proof of \eqref{lemma:chain:map}:}
Using \eqref{left:action}, we have 
\begin{align*}
 {\XO} \cdot 
\sum_{i} \AxA{i}{1} \mywedge \cdots\mywedge \AxA{i}{m}
& =
\sum_{i}
\sum_{k=1}^{m}
 \AxA{i}{1} \mywedge \cdots\mywedge \Sbt {\XO}{\AxA{i}{k}}
\mywedge  \cdots \mywedge  \AxA{i}{m} \\
&= \sum_{i}
\sum_{k=1}^{m}
 \AxA{i}{1} \mywedge \cdots\mywedge (- w( \AxA{i}{k} ) + 
  \bar{h}( \AxA{i}{k} )  ){\AxA{i}{k}}
\mywedge  \cdots \mywedge  \AxA{i}{m} \\
&= \sum_{i}
\sum_{k=1}^{m} ( - w( \AxA{i}{k} ) + \bar{h}( \AxA{i}{k} )  )
 \AxA{i}{1} \mywedge \cdots\mywedge 
 {\AxA{i}{k}}
\mywedge  \cdots \mywedge  \AxA{i}{m} \\
&= \sum_{i}
 ( -w + h ) \AxA{i}{1} \mywedge \cdots\mywedge 
 {\AxA{i}{k}}
\mywedge  \cdots \mywedge  \AxA{i}{m} \\
&= ( -w + h ) 
\sum_{i} \AxA{i}{1} \mywedge \cdots\mywedge 
 {\AxA{i}{k}}
\mywedge  \cdots \mywedge  \AxA{i}{m} \;. 
\end{align*}
\kmqed

Using Lemma above, we have the next proposition: 
\begin{prop}
Define a map \(\ds\phi : \myCSW{m}{w,h} \to \myCSW{m+1}{w,h}\) by
\(\phi(U)= \XO \mywedge U\). Then we have
\begin{equation}
\pdel \circ \phi + \phi\circ \pdel = (-w+h) \mathop{id}\;. 
\end{equation}
\end{prop} 
\textbf{Proof:} Take \(\ds \forall W \in \myCSW{m}{w,h}\). 
\begin{align*}
\pdel ( \phi W ) =& \pdel( \XO\mywedge W ) 
\mathop{ = }^{\eqref{pdel:recurs:left}} - \XO \mywedge \pdel W + \XO\cdot W \mathop{=}^{\eqref{lemma:chain:map}}
- \phi( \pdel W ) + (-w+h) W
\end{align*} \kmqed

Directly from this proposition we have the next theorem.  
\begin{thm} [m-th Betti number]\label{thm:m:th:Betti}
Each $m$-th Betti numbers
of \((w,h)\)-weighted chain complex \(\ds\{ \myCSW{\bullet}{w,h} \}\) 
is 0 if \(w \ne h\). 
\end{thm}

\textbf{Proof:} 
Take a general cycle 
\(\ds W \in \myCSW{m}{w,h}\).  
Th proposition above yields
\[ (-w+h) W = \pdel ( \phi(W) ) + \phi(\pdel U) = \pdel ( \phi(W)
)\quad\text{and}\quad 
W 
 = \frac{1 }{-w+h} \pdel ( \XO \mywedge W)\quad\text{if}\quad w\ne
 h\;.\] 
 \kmqed

\begin{remark} When $w=h$,  
Theorem \ref{thm:m:th:Betti} says \( \XO\mywedge U\) is a  cycle 
if \(U\) is a cycle in \(\ds \myCSW{m}{w,h}\). 
\end{remark}

\begin{remark} 
We have the table of Betti numbers of  
\(\ds \{ \myCSW{\bullet}{0,0} \}\) of \(\ds\mR^{2}\)
in Example \ref{exam:n2:w0}, 
which shows non-trivial Betti numbers:  
\( b_{0} =1, b_{5} =2, b_{7} =1, b_{8} =2\).  

\kmcomment{
\begin{center}
\begin{tabular}{c | *{9}{c} }
m & 0 & 1& 2 & 3 & 4 & 5 & 6 & 7 & 8 \\\hline
$\dim$ &1& 4 &  18 & 60 & 120 & 156 & 134  & 68 & 15 
\\
$\dim\pdel$ & 0 & 4 & 14 & 46 & 74 & 80 & 54 & 13 & 0 
\\\hline
Betti & 1 &  0 & 0 & 0 & 0 & 2 & 0 & 1 & 2 
\end{tabular}
\end{center} 
$m=0$ is possible because of $w=h=0$, and \(\ds \myCSW{0}{0,0}  = \mR\).   
The Euler number is $0$. 
}
\end{remark}

Theorem \ref{thm:m:th:Betti} being concerned with $m$-th Betti numbers
makes sense for $m=1$ with assumption \(w\ne h\). But we have the next
result without any restriction for the first Betti number.  

\begin{thm} [1st Betti number] \label{thm:1:th:Betti}
The first Betti number 
of \((w,h)\)-weighted chain complex \(\ds\{ \myCSW{\bullet}{w,h} \}\) 
is 0 for each double weight \((w,h)\). 
\end{thm}
\textbf{Proof:}
Fix general weight \((w,h)\). 
\kmcomment{
The weighted chain space \(\myCSW{m}{w}\) corresponds to \(\ds
\mynabla{w+m}{m}\), in particular, 
 \(\myCSW{1}{w} = \mynabla{w+1}{1} = [ (w+1)^{1} ]\) and 
 \(\myCSW{2}{w} = \mynabla{w+2}{2} = [ 1^{1}, (w+1)^{1} ]
\sqcup   [ 2^{1}, {w}^{1} ]
\sqcup   [ 3^{1}, (w-2)^{1} ]  
\sqcup \cdots   
 \) where the shape of end of the sequence depends on $w$ is even or odd.  
Thus,} 
\(\ds \myCSW{1}{w,h} = \mathfrakX{w+1}{h+1}\) and 
\(\ds \myCSW{2}{w,h}= \sum_{\substack{p+q = 2+w \\ a+b=2+h} }\mathfrakX{p}{a}
\mywedge  \mathfrakX{q}{b} \).  

Take \(\forall\, U \in  \mathfrakX{w+1}{h+1}\). Then 
\(\ds  \bibun{k}\mywedge ( x_{k} U )\in 
 \mathfrakX{1}{0} \mywedge \mathfrakX{w+1}{h+2} \subset 
 \myCSW{2}{w,h} \).  Now we see 
\begin{align*}
\pdel (
 \bibun{k}\mywedge x_{k} U  
) =& \Sbt{ \bibun{k} } { x_{k} U } = U + x_{k} 
 \Sbt{ \bibun{k} } { U }  \\
\sum_{k=1}^{n}\pdel (
 \bibun{k}\mywedge x_{k} U  
) =&  n U + \sum_{k=1}^{n} x_{k} \Sbt{ \bibun{k} } { U }  
=  (n + 1 + h) U \; . 
\end{align*}
\kmqed


We have result about the second Betti number when $w=h=0$. 
\begin{prop} 
The 2nd Betti number is zero when  $w=h=0$. 
\end{prop}
\kmcomment{
\textbf{Proof:}
Now, \(\myCSW{2}{0,0} =
\sum_{a+b=0+2, a\leq b} \mathfrakX{1}{a} \mywedge  \mathfrakX{1}{b} 
= \mathfrakX{1}{0} \mywedge  \mathfrakX{1}{2} 
+ 
\mathfrakX{1}{1} \mywedge  \mathfrakX{1}{1} \).  
Take \( \sum_{i}Y_{i,0} \mywedge Y_{i,2} 
\in  \mathfrakX{1}{0} \mywedge  \mathfrakX{1}{2} \) and 
\( \sum_{j} Y_{j} \mywedge Z_{j} 
\in  \mathfrakX{1}{1} \mywedge  \mathfrakX{1}{1} \) with the assumption  
\( \pdel (
\sum_{i} Y_{i,0} \mywedge Y_{i,2}  + \sum_{j}  Y_{j} \mywedge Z_{j})
=0\;,\; \text{i.e.,}\; \sum_{i}
 \Sbt{ Y_{i,0}}{ Y_{i,2} } +  \sum_{j}
 \Sbt{ Y_{j}}{ Z_{j} }  = 0 \). Now consider \( \sum_{i}
 \bibun{k} \mywedge Y_{i,0} \mywedge ( x_{k} Y_{i,2} ) 
 + \sum_{j} \bibun{k} \mywedge Y_{j} \mywedge ( x_{k} Z_{j} )\in\myCSW{3}{0,0}\). 
\begin{align*}
& \pdel (
 \sum_{i} \bibun{k} \mywedge Y_{i,0} \mywedge ( x_{k} Y_{i,2} ) 
 + \sum_{j} \bibun{k} \mywedge Y_{j} \mywedge ( x_{k} Z_{j} )
) 
\\ = & \sum_{i} \left(  - \bibun{k} \mywedge \Sbt{Y_{i,0}} { x_{k} Y_{i,2} } 
+ \Sbt{\bibun{k}}{Y_{i,0}} \mywedge ( x_{k} Y_{i,2} ) 
+ {Y_{i,0}} \mywedge 
 \Sbt{\bibun{k}} { x_{k} Y_{i,2} } \right) \\
& + 
\sum_{j} \left(  - \bibun{k} \mywedge \Sbt{Y_{j}} { x_{k} Z_{j} } 
+ \Sbt{\bibun{k}}{Y_{j}} \mywedge ( x_{k} Z_{j} ) 
+ {Y_{j}} \mywedge \Sbt{\bibun{k}} { x_{k} Z_{j} } \right) \\
= & \sum_{i} \left(  - \bibun{k} \mywedge 
( \Sbt{Y_{i,0}} { x_{k}}  Y_{i,2} + x_{k} \Sbt{Y_{i,0}} {  Y_{i,2} } )
+  \Sbt{\bibun{k}}{Y_{i,0}} \mywedge ( x_{k} Y_{i,2} ) 
+ {Y_{i,0}} \mywedge 
 \Sbt{\bibun{k}} { x_{k} Y_{i,2} } \right) \\
& + 
\sum_{j} \left(  - \bibun{k} \mywedge 
( \Sbt{Y_{j}} { x_{k} } { Z_{j} } + x_{k} \Sbt{Y_{j}} {  Z_{j} } )
+ \Sbt{\bibun{k}}{Y_{j}} \mywedge ( x_{k} Z_{j} ) 
+ {Y_{j}} \mywedge \Sbt{\bibun{k}} { x_{k} Z_{j} } \right) \\ 
\noalign{from cycle condition, we have}
= & \sum_{i} \left(  - \bibun{k} \mywedge 
( \Sbt{Y_{i,0}} { x_{k}}  Y_{i,2}  )
+  \Sbt{\bibun{k}}{Y_{i,0}} \mywedge ( x_{k} Y_{i,2} ) 
+ {Y_{i,0}} \mywedge 
 \Sbt{\bibun{k}} { x_{k} Y_{i,2} } \right) \\
& + 
\sum_{j} \left(  - \bibun{k} \mywedge 
( \Sbt{Y_{j}} { x_{k} } { Z_{j} }  )
+ \Sbt{\bibun{k}}{Y_{j}} \mywedge ( x_{k} Z_{j} ) 
+ {Y_{j}} \mywedge \Sbt{\bibun{k}} { x_{k} Z_{j} } \right)  \\ 
\noalign{since \(\ds 
 \Sbt{Y_{i,0}} { x_{k}}  \) are constant number, we have} 
= & \sum_{i} \left(  - \Sbt{Y_{i,0}}{x_{k}} \bibun{k} \mywedge Y_{i,2}  
+ 0 
+ {Y_{i,0}} \mywedge 
 \Sbt{\bibun{k}} { x_{k} Y_{i,2} } \right) \\
& + 
\sum_{j} \left(  - \bibun{k} \mywedge 
( \Sbt{Y_{j}} { x_{k} } { Z_{j} }  )
+ \Sbt{\bibun{k}}{Y_{j}} \mywedge ( x_{k} Z_{j} ) 
+ {Y_{j}} \mywedge \Sbt{\bibun{k}} { x_{k} Z_{j} } \right) \;. \\ 
\noalign{Now, summing up by $k$, we have}
& \sum_{k} \pdel \left(
 \sum_{i} \bibun{k} \mywedge Y_{i,0} \mywedge ( x_{k} Y_{i,2} ) 
 + \sum_{j} \bibun{k} \mywedge Y_{j} \mywedge ( x_{k} Z_{j} )
\right) \\
= & - \sum_{i} 
 Y_{i,0} \mywedge Y_{i,2}   
+ \sum_{i,k} {Y_{i,0}} \mywedge (
 {  Y_{i,2} } 
 + x_{k} \Sbt{\bibun{k}} {  Y_{i,2} } 
 ) \\
& + 
\sum_{j,k} \left(  - \bibun{k} \mywedge 
( \Sbt{Y_{j}} { x_{k} } { Z_{j} }  )
+ \Sbt{\bibun{k}}{Y_{j}} \mywedge ( x_{k} Z_{j} ) 
+ {Y_{j}} \mywedge ( {  Z_{j} } 
+ x_{k} \Sbt{\bibun{k}} {  Z_{j} } 
)\right)\;  \\ 
= & - \sum_{i}  {Y_{i,0}} \mywedge Y_{i,2}   
+ \sum_{i} {Y_{i,0}} \mywedge (
 (n+2) {  Y_{i,2} } 
 ) \\
&  
- \sum_{j,k} \bibun{k} \mywedge 
( \Sbt{Y_{j}} { x_{k} } { Z_{j} }  )
+ \sum_{j,k}\Sbt{\bibun{k}}{Y_{j}} \mywedge ( x_{k} Z_{j} ) 
+ \sum_{j}{Y_{j}} \mywedge ( (n+1) {  Z_{j} } 
)\;  \\ 
= & 
 (n+1)(   \sum_{i} {Y_{i,0}} \mywedge   {  Y_{i,2} } 
+ \sum_{j}{Y_{j}} \mywedge   {  Z_{j} } )
   - \sum_{j,k} \bibun{k} \mywedge 
( \Sbt{Y_{j}} { x_{k} } { Z_{j} }  )
+ \sum_{j,k}\Sbt{\bibun{k}}{Y_{j}} \mywedge ( x_{k} Z_{j} ) 
\;.    
\end{align*}
We show the sum of the last two terms is zero as follows: 
Since \(\ds Y_{j}\in \mathfrakX{1}{1}\),  \( Y_{j} = \sum_{k} Y^{k}_{j}
\bibun{k}
= \sum_{k,\ell} c^{k,\ell}_{j} x_{\ell} \bibun{k}
\) where \(\ds c^{k,\ell}_{j}\) are constant. 
Since 
\begin{align*}  \Sbt{Y_{j}}{x_{k}} Z_{j} =& 
Y^{k}_{j} Z_{j} = 
\sum_{\ell} c^{k,\ell}_{j} x_{\ell} Z_{j}\\
\noalign{and}
\Sbt{\bibun{k}}{Y_{j}} =& \sum_{\ell} \Sbt{\bibun{k}}{ Y^{\ell}_{j} }
\bibun{\ell}
= \sum_{\ell} c^{\ell,k}_{j} \bibun{\ell} \;,
\\\noalign{we get} 
\text{2nd term} + \text{3rd term} =& - \sum_{j,k} \bibun{k} \mywedge 
(  
\sum_{\ell} c^{k,\ell}_{j} x_{\ell} Z_{j}
)
+ \sum_{j,k} 
 \sum_{\ell} c^{\ell,k}_{j} \bibun{\ell} \mywedge ( x_{k} Z_{j} ) = 0\;.
\end{align*} 
\kmqed
}
\textbf{Proof:} 
Since \(\myCSW{2}{0,0} =
\sum_{a+b=0+2, a\leq b} \mathfrakX{1}{a} \mywedge  \mathfrakX{1}{b} 
= \mathfrakX{1}{0} \mywedge  \mathfrakX{1}{2} 
+ \mathfrakX{1}{1} \mywedge  \mathfrakX{1}{1} \),  
a general element \( T \in \myCSW{2}{0,0} \)
is given by 
\( T = \sum_{i} A_{i} \mywedge B_{i} + \sum_{j<\ell}p^{j,\ell} Y_{j} \mywedge Y_{\ell} \) where
\[ A_{i}\in \mathfrakX{1}{0}\;,\;  
 B_{i}\in \mathfrakX{1}{2}\;,\;  
  Y_{j} \in \mathfrakX{1}{1}\;,\;
   p^{j,\ell} + p^{\ell,j} = 0  \] 
and 
 satisfies cycle condition 
\( \pdel 
(\sum_{i}A_{i}\mywedge B_{i}+\sum_{j<\ell}p^{j,\ell}Y_{j}\mywedge Y_{\ell})
=0\), i.e., \(\sum_{i}
 \Sbt{ A_{i}}{ B_{i} } +  \sum_{j<\ell}p^{j,\ell}
 \Sbt{ Y_{j}}{ Y_{\ell} }  = 0 \). 
 
Consider \[ \sum_{i}
 \bibun{k} \mywedge A_{i} \mywedge ( x_{k} B_{i} ) 
 + \sum_{j<\ell} p^{j,\ell} \bibun{k} \mywedge Y_{j} \mywedge ( x_{k}
 Y_{\ell} )\in\myCSW{3}{0,0}\;. \] 
\begin{align*} 
& \pdel \left(
 \sum_{i} \bibun{k} \mywedge A_{i} \mywedge ( x_{k} B_{i} ) 
 + \sum_{j<\ell}p^{j,\ell} \bibun{k} \mywedge Y_{j} \mywedge ( x_{k} Y_{\ell} )
\right) 
\\ = &  - \sum_{i} \bibun{k} \mywedge \Sbt{A_{i}} { x_{k} B_{i} } 
+ \sum_{i} \Sbt{\bibun{k}}{A_{i}} \mywedge ( x_{k} B_{i} ) 
+ \sum_{i} {A_{i}} \mywedge 
 \Sbt{\bibun{k}} { x_{k} B_{i} } 
 \\ &  
-\sum_{j<\ell}p^{j,\ell}  \bibun{k} \mywedge \Sbt{Y_{j}} { x_{k} Y_{\ell} } 
+ \sum_{j<\ell}p^{j,\ell} \Sbt{\bibun{k}}{Y_{j}} \mywedge ( x_{k} Y_{\ell} ) 
+ \sum_{j<\ell}p^{j,\ell}Y_{j}\mywedge  \Sbt{\bibun{k}}{x_{k}Y_{\ell}}  
 \\
= &  
  - \sum_{i}\bibun{k} \mywedge 
( \Sbt{A_{i}} { x_{k}}  B_{i} + x_{k} \Sbt{A_{i}} {B_{i}} )
+ \sum_{i} \Sbt{\bibun{k}}{A_{i}} \mywedge ( x_{k} B_{i} ) 
+ \sum_{i}{A_{i}} \mywedge \Sbt{\bibun{k}} { x_{k} A_{i} }  \\
& 
  - \sum_{j<\ell} p^{j,\ell} 
  \bibun{k} \mywedge 
( \Sbt{Y_{j}} { x_{k} } { Y_{\ell} } + x_{k} \Sbt{Y_{j}} {  Y_{\ell} } )
+\sum_{j<\ell}  p^{j,\ell}  \Sbt{\bibun{k}}{Y_{j}} \mywedge ( x_{k} Y_{\ell} ) 
+\sum_{j<\ell}  p^{j,\ell}  {Y_{j}} \mywedge \Sbt{\bibun{k}} { x_{k} Y_{\ell} }  \\ \noalign{from cycle condition, we have}
= & 
- \sum_{i}\bibun{k} \mywedge ( \Sbt{A_{i}} { x_{k}}  B_{i}  )
+ \sum_{i} \Sbt{\bibun{k}}{A_{i}} \mywedge ( x_{k} B_{i} ) 
+ \sum_{i}{A_{i}} \mywedge \Sbt{\bibun{k}} { x_{k} B_{i} }  \\
& - \sum_{j<\ell} p^{j,\ell}  
\bibun{k} \mywedge 
( \Sbt{Y_{j}} { x_{k} } { Y_{\ell} }  )
+ \sum_{j<\ell} p^{j,\ell}  \Sbt{\bibun{k}}{Y_{j}} \mywedge ( x_{k} Y_{\ell} ) 
+ \sum_{j<\ell} p^{j,\ell}  {Y_{j}} \mywedge \Sbt{\bibun{k}} { x_{k} Y_{\ell} }  \\ 
\noalign{since \(\ds 
 \Sbt{A_{i}} { x_{k}}  \) are constant number, we have} 
= & 
- \sum_{i}\Sbt{A_{i}} { x_{k}}\bibun{k} \mywedge    B_{i}  
+ 0 
+ \sum_{i}{A_{i}} \mywedge \Sbt{\bibun{k}} { x_{k} B_{i} }  \\
& - \sum_{j<\ell} p^{j,\ell}  
\bibun{k} \mywedge 
( \Sbt{Y_{j}} { x_{k} } { Y_{\ell} }  )
+ \sum_{j<\ell} p^{j,\ell}  \Sbt{\bibun{k}}{Y_{j}} \mywedge ( x_{k} Y_{\ell} ) 
+ \sum_{j<\ell} p^{j,\ell}  {Y_{j}} \mywedge \Sbt{\bibun{k}} { x_{k}
Y_{\ell} }\;.  \\ 
\noalign{Now, summing up by $k$, we have}
& \sum_{k} \pdel \left(
 \sum_{i} \bibun{k} \mywedge A_{i} \mywedge ( x_{k} B_{i} ) 
 + \sum_{j} \bibun{k} \mywedge Y_{j} \mywedge ( x_{k} Y_{\ell} )
\right) \\
= & - \sum_{i} 
 A_{i} \mywedge B_{i}   
+ \sum_{i,k} {A_{i}} \mywedge (
 {  B_{i} } 
 + x_{k} \Sbt{\bibun{k}} {  B_{i} } 
 ) \\
& + 
\sum_{j,k} \left(  - \bibun{k} \mywedge 
( \Sbt{Y_{j}} { x_{k} } { Y_{\ell} }  )
+ \Sbt{\bibun{k}}{Y_{j}} \mywedge ( x_{k} Y_{\ell} ) 
+ {Y_{j}} \mywedge ( {  Y_{\ell} } 
+ x_{k} \Sbt{\bibun{k}} {  Y_{\ell} } 
)\right)\;  \\ 
= & - \sum_{i}  {A_{i}} \mywedge B_{i}   
+ \sum_{i} {A_{i}} \mywedge ( (n+2) {  B_{i} } ) \\
&  - \sum_{j<\ell,k} p^{j,\ell}\bibun{k} \mywedge 
( \Sbt{Y_{j}} { x_{k} } { Y_{\ell} }  )
+ \sum_{j<\ell,k}p^{j,\ell}\Sbt{\bibun{k}}{Y_{j}} \mywedge ( x_{k} Y_{\ell} ) 
+ \sum_{j<\ell} p^{j,\ell}{Y_{j}} \mywedge ( (n+1) {  Y_{\ell} } 
)\;  \\ 
= & 
 (n+1)(   \sum_{i} {A_{i}} \mywedge   {  B_{i} } 
+ \sum_{j<\ell} p^{j,\ell} {Y_{j}} \mywedge   {  Y_{\ell} } )
\\& \qquad 
   - \sum_{j<\ell,k}p^{j,\ell} \bibun{k} \mywedge 
( \Sbt{Y_{j}} { x_{k} } { Y_{\ell} }  )
+ \sum_{j<\ell,k}p^{j,\ell}\Sbt{\bibun{k}}{Y_{j}} \mywedge ( x_{k} Y_{\ell} ) 
\;.    
\end{align*}
We show the sum of the last two terms is zero as follows: 
Since \(\ds  \mathfrakX{1}{1}\ni  Y_{j}  
= \sum_{k,\ell} Y^{k,\ell}_{j} x_{\ell} \bibun{k}
\) where \(\ds Y^{k,\ell}_{j}\) are constant.  
\begin{align*}  \Sbt{Y_{j}}{x_{k}} Y_{\ell} =& 
Y^{k}_{j} Y_{\ell} = 
\sum_{t} Y^{k,t}_{j} x_{t} Y_{\ell}\;,\quad 
\Sbt{\bibun{k}}{Y_{j}} = \sum_{s} \Sbt{\bibun{k}}{ Y^{s}_{j} }
\bibun{s}
= \sum_{s} Y^{s,k}_{j} \bibun{s} 
\\\noalign{we get 
\text{2nd term} + \text{3rd term} is }
& - \sum_{j<\ell,k}p^{j,\ell} 
\bibun{k} \mywedge (  \sum_{t} Y^{k,t}_{j} x_{t} Y_{\ell})
+ \sum_{j<\ell,k}p^{j,\ell} 
 \sum_{s} Y^{s,k}_{j} \bibun{s} \mywedge ( x_{k} Y_{\ell} )\\
=&
- \sum_{j<\ell}p^{j,\ell} \sum_{k}
\bibun{k} \mywedge (  \sum_{t} Y^{k,\ell}_{j} x_{\ell} Y_{\ell})
+ \sum_{j<\ell}p^{j,\ell} 
 \sum_{s} Y^{s,k}_{j} \bibun{s} \mywedge ( x_{k} Y_{\ell} )\\
=&
- \sum_{j<\ell}p^{j,\ell} \sum_{k}
\bibun{k} \mywedge (  \sum_{t} Y^{k,\ell}_{j} x_{\ell} Y_{\ell})
+ \sum_{j<\ell}p^{j,\ell} 
 \sum_{s} 
 \bibun{s} \mywedge ( Y^{s,k}_{j} x_{k} Y_{\ell} )\\
 =& 
 0\;.
\end{align*} 
\kmqed

\begin{remark}
We expect to know the second Betti number of the chain complex \(\ds \{\myCSW{\bullet}{w,w}\}\)
for \(w>0\). 
We know the second Betti number is 0  for lower $n$.   
\end{remark}

\kmcomment{
\section{Cohomology groups of pre Lie superalgebra}
In this section, we only recall the definition of cohomology theory for
pre Lie superalgebra \(\ds \frakg = \sum_{i} \frakg_{i}\) 
 with \(\frakg\)-module $V$.

If a multilinear map
\( \sigma : \underbrace{\frakg \otimes \cdots \otimes
\frakg}_{m-\text{times}}  \rightarrow V^{*} \) 
satisfies the following ``super skew-symmetric property'':
\[
\sigma(\cdots \otimes Y_{k}\otimes Y_{k+1}\otimes 
\cdots ) = (-1)^{ 
1+y_{k}y_{k+1} }
\sigma(\cdots\otimes\mathop{Y_{k+1}}_{<k>}\otimes
\mathop{Y_{k}}_{<k+1>} \otimes \cdots)\;,\] 
\kmcomment{
\[
\sigma(\cdots \otimes Y_{i}\otimes \cdots \otimes Y_{j}\otimes\cdots) = (-1)^{ 
(1+y_{i})(1+y_{j}) + (y_{i}+y_{j}) \sum_{s=i}^{j} y_{s}  }
\sigma(\cdots\otimes\mathop{Y_{j}}_{<i>}\otimes\cdots\otimes
\mathop{Y_{i}}_{<j>} \otimes \cdots)\;,\] 
}
then we call \(\sigma\) is $m$-th cochain,  
and denote the space of all $m$-th cochains by \(\ds \myCC{m}{} \).  
Define a map \(\ds \myd_{V}\)  for $m$-th cochain \(\sigma \) by 
\[ (\myd_{V} \sigma )
(Y_{1}\otimes \ldots \otimes Y_{m+1}) = 
\sigma ( \pdel ( 
Y_{1}\mywedge \cdots \mywedge   Y_{m+1}))
+ (-1)^{m} \sum_{i} (-1)^{ \sum_{s>i} (1+y_{i} y_{s}) } Y_{i}\cdot
\sigma( \CYm{i}{m+1} ) 
\;.   
\] where the last term \(\ds Y_{i}\cdot \sigma\) is the dual action of
\(\frakg\) on \(V\).  

We see that 
\( \myd_{V}  \sigma \) also 
satisfies the ``super skew-symmetric property''
and have a map
\[ \myd_{V} : \myCC{m}{ } \rightarrow \myCC{m+1}{ }\;. \]
\textbf{Proof:}
\begin{align*}
& (\pdel_{V} \sigma)(\cdots\otimes\mathop{Y_{k+1}}_{<k>}\otimes
\mathop{Y_{k}}_{<k+1>} \otimes \cdots \otimes Y_{m+1} )
\\
=& 
\sigma( \pdel( 
\cdots\otimes\mathop{Y_{k+1}}_{<k>}\otimes
\mathop{Y_{k}}_{<k+1>} \otimes \cdots \otimes Y_{m+1} ))
\\& + (-1)^{m} \sum_{\ell<k} (-1)^{ \sum_{s>\ell} (1+y_{\ell} y_{s}) } Y_{\ell}\cdot
\sigma(\dots\otimes \widehat{Y_{\ell}}\otimes \dots\otimes Y_{k+1}\otimes
Y_{k}\otimes \dots ) 
\\& + (-1)^{m} \sum_{\ell>k+1} (-1)^{ \sum_{s>k+1} (1+y_{\ell} y_{s}) } Y_{\ell}\cdot
\sigma(\dots\otimes Y_{k+1}\otimes Y_{k}\otimes \dots 
\otimes \widehat{Y_{\ell}}\otimes \dots) 
\\& + (-1)^{m} (-1)^{1+y_{k}y_{k+1} +  \sum_{s>k+1} (1+y_{k+1} y_{s}) } 
Y_{k+1}\cdot
\sigma(\dots\otimes Y_{k}\otimes 
 \widehat{Y_{k+1}}\otimes 
\dots ) 
\\& + (-1)^{m} (-1)^{1+y_{k}y_{k+1} +  \sum_{s>k} (1+y_{k} y_{s}) } 
Y_{k}\cdot
\sigma(\dots\otimes \widehat{Y_{k}} \otimes 
 Y_{k+1}\otimes \dots ) 
\\ = &  (-1)^{ 1+ y_{k+1} y_{k} } 
(\pdel_{V} \sigma)(\cdots\otimes\mathop{Y_{k}}_{<k>}\otimes
\mathop{Y_{k+1}}_{<k+1>} \otimes \cdots \otimes Y_{m+1} ) \;. 
\end{align*}
Since \(\ds \myd_{V} \) is the dual map of \(\ds \pdel_{V}\), 
we see that \(\ds \myd_{V} \circ \myd_{V} = 0\). 
Thus, we have the
cohomology group \(\ds \myHH{m}{ }(\frakg,V^{*})\)  of a super Lie
algebra \(\frakg\). Also, we restrict cochain spaces by weight or
double-weight and have those cohomology groups, too. 

We may consider the ``formal dual'' \(\ds \frakg^{*}\) of \(\frakg\),
and 
\[\ds \myCC{m}{w} = \sum_{
\substack{
i_{1} \leq \cdots \leq i_{m}\\
\sum_{s} i_{s} = w 
} }
\frakg_{i_1}^{*} \mywedge \cdots  \mywedge 
\frakg_{i_m}^{*} \otimes V^{*} \] is the $m$-th cochain space with weight
$w$. 
If the pre Lie superalgebra is double-weighted then 
\[\ds \myCC{m}{w,h} = \sum_{
\substack{
i_{1} \leq \cdots \leq i_{m}\;,\;
\sum_{s} i_{s} = w \\
\sum_{s} h_{s} = h 
} }
\frakg_{i_1,h_1}^{*} \mywedge \cdots  \mywedge 
\frakg_{i_m,h_m}^{*} \otimes V^{*} \] is the $m$-th cochain space with
double-weight \(w,h\). 
}

\nocite{MR1365257}
\nocite{GB:J:Merker}

\bibliographystyle{plain}
\bibliography{km_refs}

\end{document}